\documentclass{llncs}

\usepackage{graphicx,epsfig}
\usepackage{amsmath}
\usepackage{amssymb}
\spnewtheorem*{mynote}{Note}{\itshape}{\rmfamily}
\long\def\symbolfootnote[#1]#2{\begingroup%
\def\thefootnote{\fnsymbol{footnote}}\footnote[#1]{#2}\endgroup}

\newcommand{\bbox}{\rule{0.6em}{0.6em}}
\newcommand{\ignore}[1]{}
\newcommand{\com}[1]{}

\begin{document}
\setcounter{page}{1}

\title  {On the SIG dimension of trees under $L_{\infty}$ metric}

\author{L. Sunil Chandran\inst{1}\and Rajesh Chitnis\inst{2}\and
Ramanjit Kumar\inst{1}} \institute{Indian Institute of Science, Dept. of Computer Science and Automation, Bangalore--560012,
India.  email: \emph{\{sunil, ramanjitkumar\}@csa.iisc.ernet.in} \and Department of Computer Science, University of Maryland,
College Park, USA email: \emph{rchitnis@cs.umd.edu}} \com{Changed my address and email}

\maketitle

 {\renewcommand{\thefootnote}{}\footnote{This work was done while the second author was a summer intern at IISc}}

\begin{abstract}
Let ${\cal P}$, where $|{\cal P}| \geq 2$, be a set of points in $d$-dimensional space with a given metric $\rho$.
 For a point $p \in {\cal P}$, let $r_p$ be the distance of $p$ with respect to $\rho$
 from its nearest neighbor in ${\cal P}$.
 Let $B(p,r_p)$ be the open ball with respect to $\rho$ centered at $p$ and having the radius $r_p$.
 We define the sphere-of-influence graph $(SIG)$ of ${\cal P}$ as the intersection
 graph of the family of sets \{$B(p,r_p)\  | \ p\in {\cal P}$\}.
 Given a graph $G$, a set of points ${\cal P}_G$ in d-dimensional space with the metric
 $\rho$ is called a $d-$dimensional SIG representation of $G$, if $G$ is the SIG of ${\cal P}_G$.\\

 It is known that the absence of isolated vertices is a necessary and sufficient
 condition for a graph to have an $SIG$ representation under $L_{\infty}$ metric in some space of finite dimension.
 The $SIG$ dimension under $L_{\infty}$ metric of a graph $G$ without isolated vertices is defined to be the minimum
 positive integer $d$ such that $G$ has a $d$-dimensional $SIG$ representation under the $L_{\infty}$ metric.
 It is denoted as $SIG_\infty(G)$.\\

 We study the $SIG$ dimension of trees under $L_{\infty}$ metric and almost completely answer an open problem
 posed by Michael and Quint (Discrete Applied Mathematics: 127, pages 447-460, 2003). Let $T$ be a tree with at least two vertices. For each $v\in V(T)$,
 let leaf-degree$(v)$ denote the number of neighbors of $v$ that are leaves. We define the maximum leaf-degree as
 $\alpha(T) = \max_{x \in V(T)}$ leaf-degree$(x)$. Let $S = \{v\in V(T)\ |$\ leaf-degree$(v) = \alpha\}$. If $|S| = 1$, we
 define $\beta(T) = \alpha(T) - 1$. Otherwise define $\beta(T) = \alpha(T)$. We show that for a tree $T$,
 $SIG_\infty(T) = \lceil \log_2(\beta + 2)\rceil$ where $\beta = \beta (T)$,
 provided $\beta$ is not of the form $2^k - 1$, for some positive integer $k \geq 1$.
 If $\beta = 2^k - 1$, then $SIG_\infty (T) \in \{k, k+1\}$.
 We show that both values are possible.\\
~\\
$Keywords:$ Sphere of Influence Graphs, Trees, $L_\infty$ norm, Intersection Graphs.
\end{abstract}

\pagestyle{plain} \pagenumbering{arabic} \setcounter{page}{1}

\bibliographystyle{plain}

\section{Introduction}
\label{Introduction}

Let $x$ and $y$ be two points in $d$-dimensional space.
For any vector $z$, let $z[i]$ denote its $i^{th}$ component.
For any positive integer $n$, let the set $\{1,2,\ldots,n\}$ be denoted by $[n]$.
 The distance $\rho(x,y)$ between $x$ and $y$ with respect to the $L_k$-metric
 (where $k \geq 1$ is a positive integer) is defined to be $\rho_k(x,y) = (\sum_{j=1}^d|x[j]-y[j]|^k)^{\frac{1}{k}}$.
 Note that $\rho_2(x,y)$ corresponds to the usual notion of the Euclidean distance between the two points $x$ and $y$.
 The distance between $x$ and $y$ under $L_\infty$-metric $\rho_\infty(x,y)$ is defined as $\max\{|x[i]-y[i]|:i=1,2,\ldots,d\}$.
 In this paper for the most part we will be concerned about the distance under $L_\infty$ metric and
 therefore we will abbreviate $\rho_\infty(x,y)$ as $\rho(x,y)$ in our proofs. Also, log will always refer to logarithm to the base 2.

\subsection{Open Balls and Closed Balls: }
\label{Balls}

For a point $p$ in $d$-dimensional space and a positive real number $r$,
 the open ball $B(p,r)$ under a given metric $\rho$, is a subset of $d$-dimensional space defined by $B(p,r) = \{x\ |\ \rho(x,p) < r\}$.
 For a point $p$ in $d$-dimensional space and a positive real number $r$, the closed ball $C(p,r)$ under a given metric $\rho$,
 is a subset of $d$-dimensional space defined by $C(p,r) = \{x\ |\ \rho(x,p) \leq r\}$.
An open ball in $d$-dimensional space, with respect to $L_\infty$ metric centered at $p$ and with radius $r$ is
 actually a ``$d$-dimensional cube'' defined as the cartesian product of $d$ open intervals namely
 $(p[1] - r, p[1] + r)$, $(p[2] - r, p[2] + r),$ $\ldots$ and $(p[d] - r, p[d] + r)$.
 In notation, $B(p, r) = \prod_{i=1}^d (p[i] - r, p[i] + r)$ where $\prod$ denotes the cartesian product.

\subsection{Maximum leaf-degree of a tree}
\label{Leaf-degree} Let $T = (V, E)$ be an (unrooted) tree with $|V| \geq 2$. A vertex $x$ of $T$ is called a leaf, if
$degree(x)=1$. For a vertex $x \in V$, let $P(x) = \{y \in V\ |\ y$ is adjacent to $x$ in $T$ and $y$ is a leaf$\}$. We define
the maximum leaf degree $\alpha$ of $T$ as $\alpha(T) = \max_{x \in V(T)}|P(x)|$. For our proof it is convenient to visualize
the tree $T$ as a rooted tree. Therefore we define a special rooted tree $T'$ corresponding to $T$, by carefully selecting a
root, as follows: Let $z \in V$ be such that $|P(z)| = \alpha$. Let $z{'} \in P(z)$. Let $T{'}$ be the rooted tree obtained
from $T$, by fixing $z{'}$ as root. In a rooted tree, a vertex is called a `leaf', if it has no children. For $x \in V$, let
$L(x) = \{y \in V(T)\ |\ y $ is a child of $x$ in $T{'}$ and $y$ is a leaf$\}$. We define $\beta (T) = \max_{x \in
V(T)}|L(x)|$. The relation between $\alpha(T)$ and $\beta(T)$ is summarized below. While $\beta(T)$ has the interpretation
given above in terms of the special rooted tree $T'$, we take the following as the formal definition of $\beta(T)$.
~\\
\begin{definition}
 Let $T$ be a tree of at least $2$ vertices and $\alpha(T)$ be the maximum leaf-degree of $T$. Let $S = \{v\in V(T)\ |\ v$ is a vertex of
 maximum leaf degree, i.e. $|P(v)| = \alpha\}$. Then we define $\beta(T) = \alpha(T)$ if $|S| \geq 2$ and $\beta(T) = \alpha(T) - 1$, if $|S| = 1$.
\end{definition}
Clearly for any tree with at least 2 vertices $\alpha(T) \geq 1$. Moreover if $\alpha(T) = 1$, then $|S| \geq 2$ and therefore
$\beta(T) \geq 1$ for all trees $T$ with at least 2 vertices.

\section{$\mathbf{SIG}$-representation and $\mathbf{SIG}$ dimensions}
\label{SIG definitions}

Let ${\cal P}$, where $|{\cal P}| \geq 2$, be a set of points in $d$-dimensional space with a given metric $\rho$. For a point
$p \in {\cal P}$, let $r_p$ be the distance of $p$ from its nearest neighbor in ${\cal P}$ with respect to $\rho$. Let
$B(p,r_p)$ be the open ball with respect to $\rho$ centered at $p$ and having a radius $r_p$. We define the
sphere-of-influence graph, $SIG$, of ${\cal P}$ as the intersection graph of the family of sets \{$B(p,r_p)\  | \ p\in {\cal
P}$\} i.e the graph will have a vertex corresponding to each set and two vertices will be adjacent if and only if the
corresponding sets intersect. Given a graph $G$, a set of points ${\cal P}_G$ in d-dimensional space (with the metric $\rho$)
is called a $d$-dimensional SIG representation of $G$, if $G$ is the SIG of ${\cal P}_G$. Note that if ${\cal P}$ is a
$d$-dimensional $SIG$ representation of $G$ then we are associating to each vertex $x$ of $G$ a point $p(x)$ in ${\cal P}$.
Given a graph $G$, the minimum positive integer $d$ such that $G$ has a $d$-dimensional $SIG$ representation (with respect to
the metric $\rho$) is called
 the $SIG$ dimension of $G$ (with respect to the metric $\rho$) and is denoted by $SIG_{\rho}(G)$.\\
 It is known that the absence of isolated vertices is a necessary and sufficient
 condition for a graph to have an $SIG$ representation under $L_{\infty}$ metric in some space of finite dimension~\cite{Michael}.
 The $SIG$ dimension under $L_{\infty}$ metric of a graph $G$ without isolated vertices is defined to be the minimum
 positive integer $d$ such that $G$ has a $d$-dimensional $SIG$ representation under the $L_{\infty}$ metric.
 It is denoted as $SIG_\infty(G)$. In this paper we may sometimes abbreviate this as $SIG(G)$.

\section{Literature Survey}
\label{Survey}

Toussaint introduced sphere-of-influence graphs (SIG) to model situations in pattern recognition and computer vision in
\cite{Tous 1}, \cite{Tous 2} and \cite{Tous 3}. Graphs which can be realized as $SIG$ in the Euclidean plane are considered in
\cite{Avis}, \cite{Harary}, \cite{Jacobson} and \cite{Lipman}. $SIG$ in general metric spaces are considered in \cite{Quint
General}.

Toussaint has used the sphere-of-influence graphs under $L_2$-metric to capture low-level perceptual information in certain
dot patterns. It is argued in \cite{Michael} that sphere-of-influence graphs under the $L_{\infty}$-metric perform better for
this purpose. Also, several results regarding $SIG_{\infty}$ dimension are proved in \cite{Michael}. Bounds for the
$SIG_{\infty}$ dimension of complete multipartite graphs are considered in \cite{Boyer}.

\section{Our result}
\label{Our results}

We almost completely answer the following open problem regarding $SIG_{\infty}$ dimension of trees posed in~\cite{Michael}.\\
~\\
Problem: (given in page 458 of \cite {Michael})  Find a formula for the $SIG_\infty$ dimension of a tree (say in terms of its degree
sequence and graphical parameters).\\
~\\
In this paper we prove the following theorem :\\
\textbf{Theorem :} For any tree $T$ with at least 2 vertices, $SIG_\infty(T) = \lceil \log_2(\beta + 2)\rceil$ where $\beta =
\beta (T)$,
 provided $\beta$ is not of the form $2^k - 1$, for some positive integer $k \geq 1$.
 If $\beta = 2^k - 1$, then $SIG_\infty (T) \in \{k, k+1\}$. In this case we show that both values can be achieved.\\

\noindent {\bf Remark:} For the realization of trees as sphere of influence graphs using other metrics or closed balls, see
Michael and Quint \cite {Quint General} and Jacobson, Lipman and Mcmorris \cite {Jacobson}.

\section{Lower Bound for $\mathbf {SIG_\infty(T)}$}
\label{Lower Bound}

\begin{lemma}
For all trees $T$, $SIG_\infty(T)\geq \lceil \log_2(\beta + 1)\rceil$.
\label{Lower Lemma}
\end{lemma}
\begin{proof}
 If $|V(T)| = 2$, then it is a single edge and in this case $SIG_\infty (T) = 1$ and so the theorem is true in this case. Now let us assume that $V(T) \geq$ 3.
 Let $t = SIG_\infty(T)$.
 Consider the special rooted tree $T'$ corresponding to $T$ defined in section \ref{Leaf-degree}.
 Let $z \in V$ be such that $|L(z)| = \beta$ in the rooted tree $T'$.
 Let $L(z) = \{y_1,y_2,\ldots,y_\beta\}$.
 Consider a $SIG$ representation ${\cal P}_T$ of the tree $T$ in $t$-dimensional space under $L_\infty$ metric.
 From the definition of $SIG$ representation it is clear that each vertex $x$ has to be adjacent to every vertex $y$ such that
 $y$ is a nearest point of $x$ in ${\cal P}_T$.
 Since for each $y \in L(z)$, $z$ is the only adjacent vertex it follows that $z$ is the unique nearest point to $y$ for each $y \in L(z)$.
\begin{claim}
For each $y\in L(z),\ Vol(B(y, r_y) \cap B(z, r_z)) \geq (r_z)^t$ where $Vol$ denotes the volume.
\end{claim}
\begin{proof}
 Following notation from section 1, $y[i]$ and $z[i]$ will denote the $i$th co-ordinate of $y$ and $z$ respectively.
 As $z$ is the nearest point to $y$ we have $r_y = \rho(y,z)$.
 From this we can infer that $r_z \leq \rho(y,z) = r_y$.
 Since $\rho(y,z) = r_y$ we have by the definition of $L_\infty$ metric, $|y[i] - z[i]| \leq r_y$ for $1 \leq i \leq t$.
 Without loss of generality, we may assume that $z$ is the origin, i.e. $z[i] = 0$ for $1 \leq i \leq t$.
 This means that $|y[i]| = |y[i] - z[i]| \leq \rho(y,z) = r_y$.
 Now consider the projection of $B(y,r_y)$ and $B(z, r_z)$ on the $i^{th}$ axis.
 Clearly these projections are the open intervals $(y[i] - r_y, y[i] + r_y)$ and $(-r_z,r_z)$ respectively.
 We claim that the length of the intersection of these two intervals is at least $r_z$.
 To see this  we consider the following two cases:
 \begin{itemize}
 \item If $y[i] \leq 0$, then $|y[i]| \leq r_y$ implies $0 \leq y[i] + r_y $. Since $r_z \leq r_y$ the interval $(-r_z,
     0)$ is contained in the interval $(y[i] - r_y, y[i] + r_y)$.
 \item If $y[i] > 0$ then $|y[i]|\leq r_y$ implies $y[i] - r_y \leq 0$. Since $r_z \leq r_y$ the interval $(0, r_z)$ is
     contained in the interval $(y[i] - r_y, y[i] + r_y)$.
\end{itemize}
 It follows that in both cases the length of $(y[i] - r_y, y[i] + r_y) \cap (-r_z,+r_z)$ is at least $r_z$. Also $B(y,
 r_y) \cap B(z, r_z) = \prod_{i = 1}^t[(y[i] - r_y, y[i] + r_y) \cap (-r_z, r_z)]$ where $\prod$ stand for the cartesian
 product. Therefore $Vol(B(y,r_y) \cap B(z, r_z)) \geq (r_z)^t$.\hfill$\bbox$
 \end{proof}

Let $\overline z$ be the parent of $z$ in $T'$.
 Note that $z$ always has a parent in $T'$. This is because $z$ cannot be the root of $T'$, since by the rule for selection of root for $T'$, the root
 can have only one child and this child cannot be a leaf since $V(T)\geq $ 3.
 Note that $\{\overline z, y_1, y_2, \ldots, y_\beta\}$ is an independent set in $T$ and
 therefore for $y, y' \in \{\overline z, y_1, y_2, \ldots, y_\beta\}$, $B(y, r_y) \cap B(y', r_{y'}) = \emptyset$.
 Now, noting that $Vol(B(z, r_z)) = (2r_z)^t$ we get
 $$(2r_z)^t = Vol(B(z,r_z)) \geq Vol(B(\overline z, r_{\overline z}) \cap B(z,r_z)) + \sum_{y \in L(z)} Vol(B(y,r_y) \cap B(z,r_z)) $$
 Now using the claim and noting that $|L(z)| = \beta$, we get
 $2^t(r_z)^t \geq \beta (r_z)^t + Vol(B(\overline z, r_{\overline z}) \cap B(z,r_z))$.
 Since $B(z, r_z)$ and $B(\overline z, r_{\overline z})$ are open balls and $\overline z$ is adjacent
 to $z$, $Vol(B(z, r_z) \cap B(\overline z, r_{\overline z})) > 0$. We infer that $2^t > \beta$ and therefore $2^t \geq (\beta
+ 1)$ since $\beta$ and $2^t$ are both integers. Hence $t \geq \lceil \log(\beta + 1) \rceil$.\hfill\qed
\end{proof}

\section{Upper Bound for SIG dimension of trees under $L_{\infty}$ metric}
\label{Upper Bound}

\subsection{Basic Notation}
\label{Basic Notation}

 For a non-leaf vertex $x$ of the rooted tree $T'$, let $C(x)$ denote the children of $x$. Let $L(x) = \{y \in C(x)\ |\ y$ is a leaf in $T'\}$.
 If $L(x) \neq \emptyset$, then we define $A(x) = C(x) - L(x)$. If $L(x) = \emptyset$, then select a vertex `$y_x$' from $C(x)$ and call
 it a ``pseudo-leaf'' of $x$.
 In this case we define $A(x) = C(x) - \{y_x\}$. We call elements of $A(x)$ as ``normal'' children of $x$.

\begin{figure}[!h]
\centering
\includegraphics[width=5in]{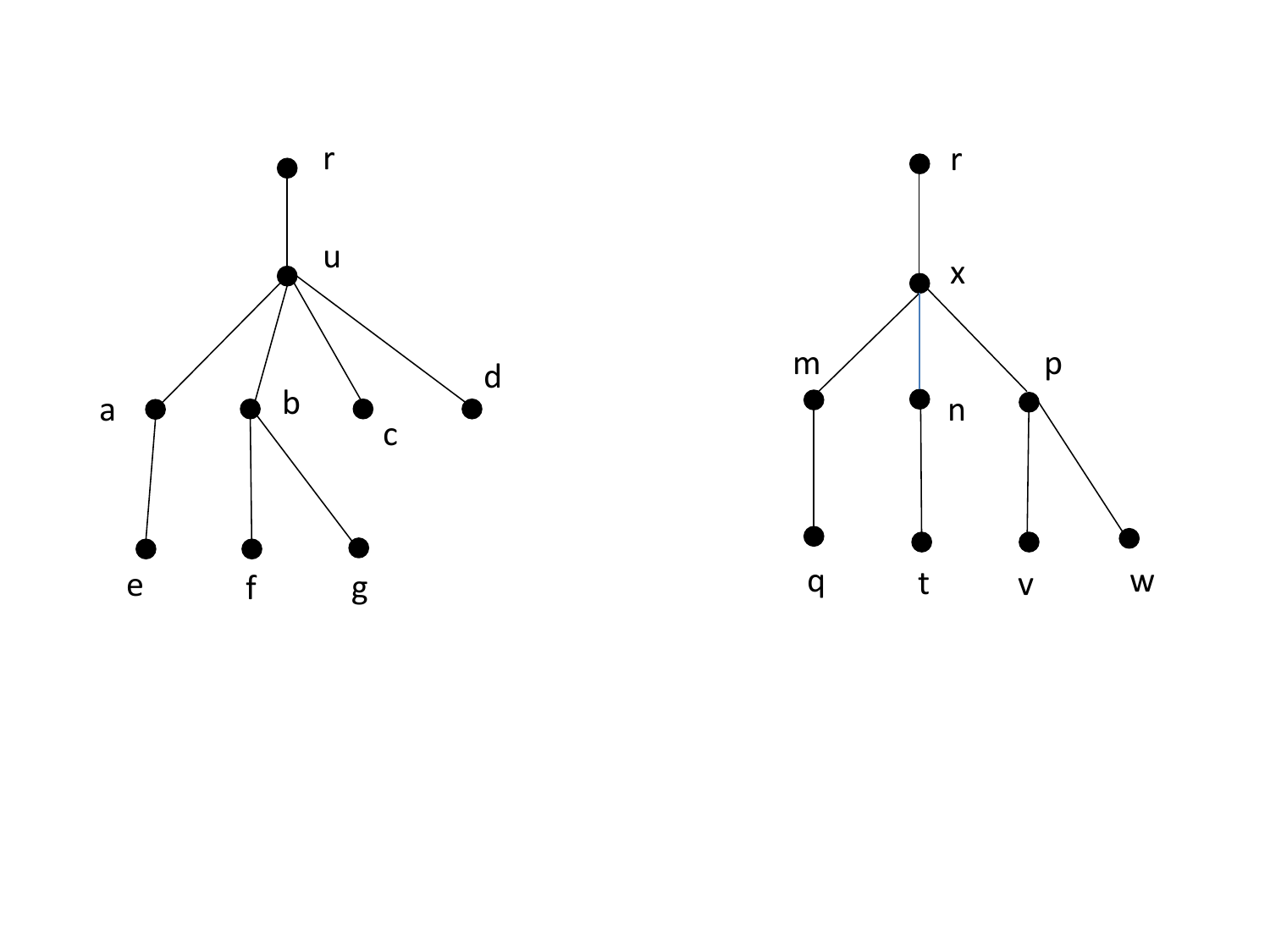}
\vspace{-30mm}
\caption
{
In the first figure $L(u)=\{c,d\}$ and $A(u)=\{a,b\}$. In the second figure $L(x)= \emptyset$ and hence we select a vertex
say $y_x = m$ as a ``pseudo-leaf" of $x$. Thus, in this case $A(x) = \{ n,p\}$.
\label{fig-normal-pseudoleaf}
}
\end{figure}

\subsection{Some More Notation Under $L_\infty$ metric}
\label{More Notation} \vspace{3mm}

\textbf{Edges and Corners of $B(p,r)$:} \\
Let ${\cal S}= \{-1,+1\}^d$ be the set of all $d$-dimensional vectors with each
component being either -1 or +1.
 The set $K(p,r)=\{p+r.S\ |\ S\in {\cal S}\}$ is the set of \emph{corners} of $B(p,r)$. Thus if $q\in K(p,r)$ then $\forall \ i\in [d]$ we
 have $q[i] = p[i] + r$ or $q[i] = p[i] - r$. Let $q$ and $q'$ be two corners of $B(p,r)$ which differ in exactly one
co-ordinate position.
 A line segment between two such corners of $B(p,r)$ is said to be an \emph{edge}\footnote
 {Note that the word {\it edge} is used in two senses in this paper: as the edge of a tree and as the edge of the ball  $B(p,r)$,
  which happens to be a $d$-dimensional cube for some $d$. Whenever
  we use this word in this paper, its meaning would be clear from the context. Moreover, in the whole of Section 6, we have taken pains
 not to use the word `edge' to indicate the tree edge, in order to avoid any possible confusion.}
 of $B(p,r)$.
 Let $q$ and $q'$ be two corners of $B(p,r)$ such that the line segment between them defines an edge. Let $i$ be such that $q[i] \neq q'[i]$.
 Then we have the following equality as sets $\{q[i], q'[i]\} = \{ p[i] - r, p[i] + r\}$.
 Note that each corner $q$, belongs to exactly $d$ edges of the ball $B(p,r)$. We denote these edges as $E_1(q),E_2(q),\ldots,E_d(q)$,
 where $E_i(q)$ is the line segment
 between $q$ and the corner $q^i$ which differs  from $q$ only in the $i$th co-ordinate.

 \vspace{3mm}

\begin{figure}[!h]
\centering
\includegraphics[width=5in]{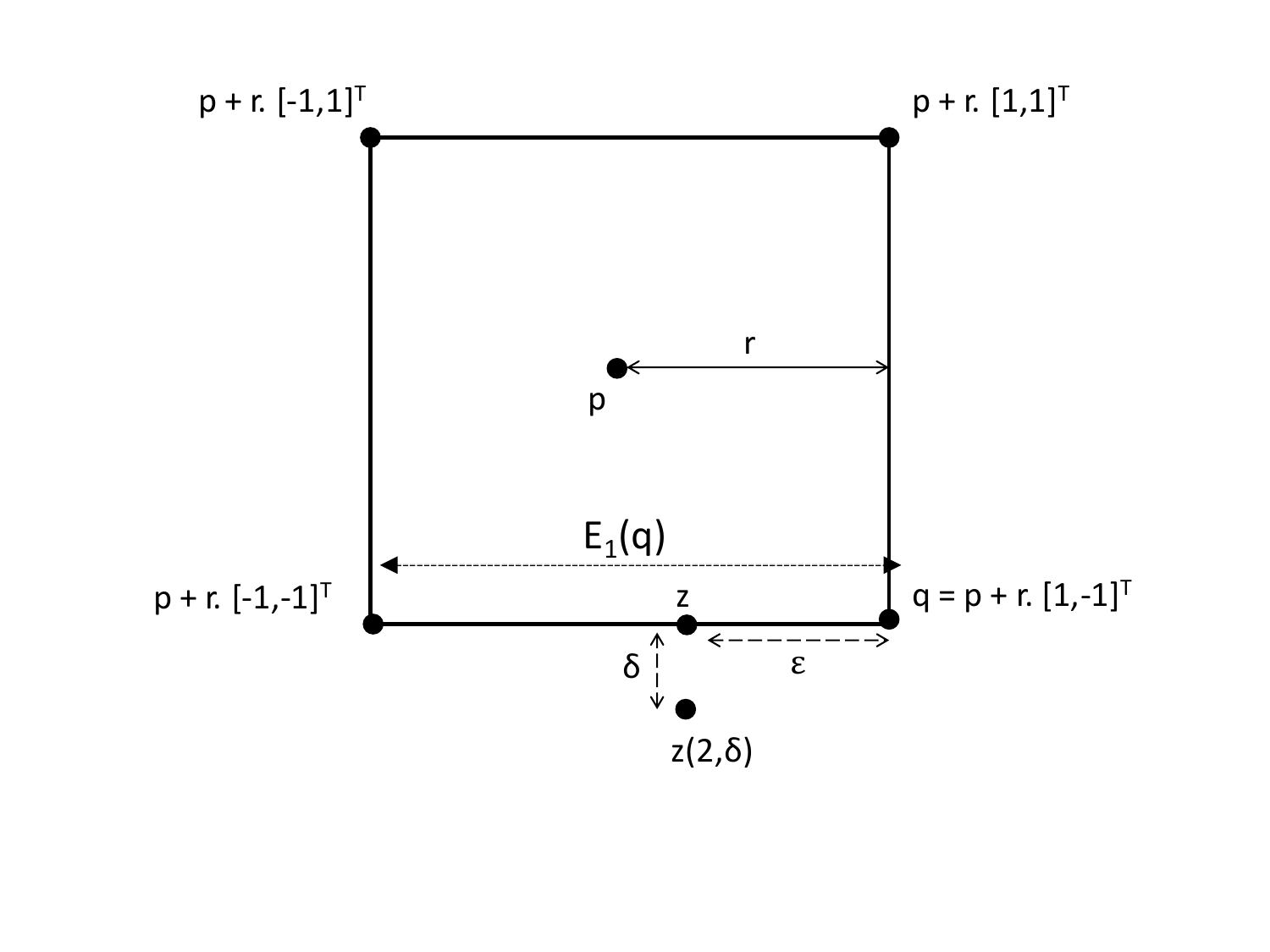}
\vspace{-10mm}
\caption
{
We give an example of the ``shifting" operation in two dimensions. We take $i=1, j=2,$ and $q$ to be the bottom right corner
of the ball $B(p,r)$. Thus, $q = p + r.S$, where $S= [1,-1]^T$.
 The point $z$ is marked on the edge $E_1(q)$.
 The point $z(j,\delta)$ is obtained by shifting $z$ along the
$j^{th}$ axis by a distance of $\delta$. Note that the first co-ordinate
of $z(j, \delta)$ is given by $z(j,\delta)[1] = z[1] = q[1] - \epsilon =
p[1] + r - \epsilon$. Also, $z(j,\delta)[2] =
z[2]+ S[2].\delta=  z[2] - \delta = q[2] - \delta
= p[2] - r - \delta$.
\label{fig-shifting}
}
\end{figure}

\noindent \textbf{The Shifting Operation:} \\
Let $q$ be a corner of $B(p,r)$ so that $q = p + r.S$ for some $S\in {\cal S}$. Consider
the edge $E_{i}(q)$.
 Let $z$ be a point on $E_i(q)$. Note that for $j\neq i$ ($1\leq j\leq d$) $z[j]=q[j]$.
 Let $\{e_1,e_2,\ldots,e_d\}$ be the canonical basis of $\mathbb{R}^d$. Now for $j\neq i$ and $\delta>0$ define $z(j,\delta)$ to be the point
 $z + (S[j].\delta)e_j$. We say that the point $z(j,\delta)$ is obtained by \emph{shifting} $z$ along the $j^{th}$ axis by the distance $\delta$.
(See figure \ref {fig-shifting}.)

\noindent \textbf{Crossing Edges:} \\
A point $p$ is said to be inside a open ball $B$ if $p\in B$, otherwise it is said to be outside
$B$. Let $B_1$ and $B_2$ be two open balls such that $B_1\cap B_2 \neq \emptyset$.
 An edge of $B_1$ is said to be \emph{crossing} with respect to the ball $B_2$ if one of the endpoint of this edge is inside $B_2$ and the other is outside.

\subsection{Algorithm to Assign Radius and Position Vector to each Vertex}
\label{Algorithm}

Let $d = \lceil \log_2(\beta+2)\rceil$ (Note that since $\beta \geq 1$, $d\geq 2$). The following algorithm assigns to each
vertex $v$ of $T'$ a point $p(v)$ (in $d$-dimensional space) called the position of $v$ and a positive real number $r(v)$
called the \textit{radius} of $v$. We will also associate with $v$ another positive real number $R(v)$ called the
\textit{super-radius} of $v$ by the following rule: If $v$ is a leaf then $R(v) = r(v)$, else $R(v) = 2r(v)$. The open ball
$B(p(v),r(v))$ will be named the \textit{ball} associated with $v$ and will be denoted by $B(v)$. The open ball $B(p(v),
R(v))$ will be called the \textit{super-ball} associated with $v$ and will be denoted by $S(v)$. Note that the
concepts `super-ball' and `super-radius' will be used only in the proofs, and therefore there will not be any
explicit mention of $S(v)$ or $R(v)$ in the algorithm.

\begin{figure}[!h]
\centering
\vspace{-7mm}
\includegraphics[width=5in]{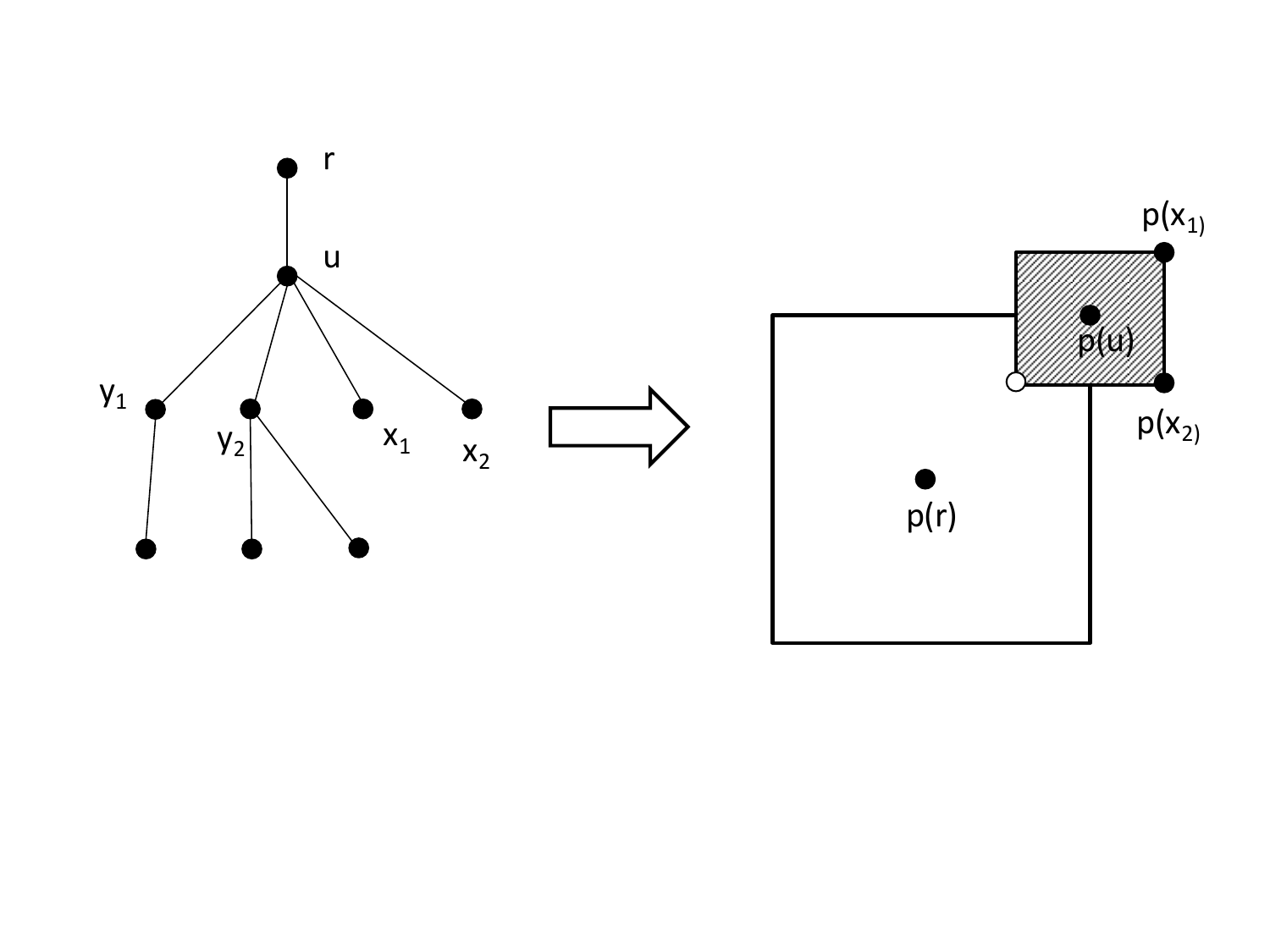}
\vspace{-25mm}
\caption
{
We illustrate Step 3.2.1 of the algorithm with a simple example in two dimensions. Note that $L(u) = \{x_1,x_2\} = C(u) - A(u)$. The
shaded box corresponds to $B(u)=B(p(u),r(u))$. The set of corners of $B(u)$
 which lie outside $B(r)$, where $r$ is the parent   of $u$,
is denoted by $K'(u)$, i.e., $K'(u) = K(u) - B(r)$. We assign different points from $K'(u)$ to the different members of $L(u)$. Thus
$p(x_1)$ and $p(x_2)$ are assigned the darkened corners of
the shaded box.  The corner of $B(u)$ inside
$B(r)$ (which is shown by the unfilled circle) is not assigned to
any member of $L(u)$.
Lemma~\ref{One Corner} states that if a vertex is a leaf or a pseudo-leaf
then exactly one corner of its box lies inside the box of its parent.
For example, in this figure the shaded box has only one corner inside
the unshaded box.
\label{first-step}
}
\end{figure}

\begin{figure}[!h]
\centering
\vspace{-7mm}
\includegraphics[width=5in]{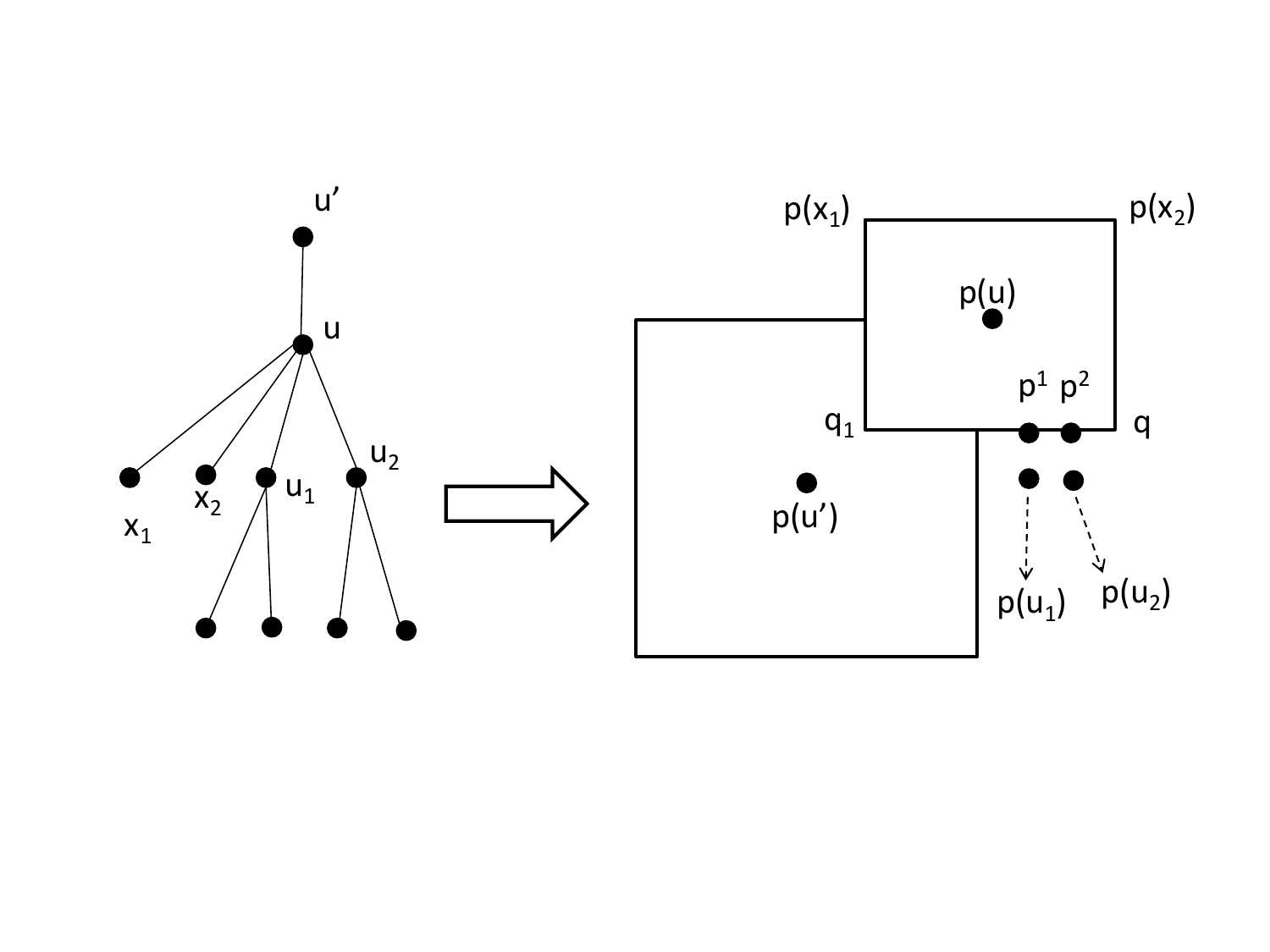}
\vspace{-25mm}
\caption
{
We illustrate the first subcase of
 Step 3.2.2 of the algorithm with a simple example
in two dimensions. For any element in $A(u)$ we first find an unused corner
from $K'(u)$. In this case
$q$ is the corner selected. Note that, when selecting the corner $q$,
 we avoided the
corners which are already assigned to $p(x_1)$ and $p(x_2)$.  Then move along the
 edge $E_1(q)$ incident to this corner to find $p^1$ and $p^2$.
Now  shift  along the $2$nd  axis to get $p(u_1$ and $p(u_2)$.
See figure \ref {second-step-part-b} to see how $B(u_i)$, for $1 \le i \le 2$
 will be placed.
\label{second-step}
}
\end{figure}

Note that the intention of the algorithm is only to assign a position vector $p(x)$ and a radius $r(x)$ to each vertex of the
tree. That the tree $T$ is the SIG of the family $\{ p(x) : x \in V(G) \}$ will be proved later. The rule to assign $r(x)$ to
each vertex $x$ is quite simple, and is given in  step 3.1 of the algorithm. The rule to assign $p(x)$ to vertex $x$ when $x$
is a leaf or a pseudo leaf of its parent is again simple, and is given in step 3.2.1 of the algorithm. The rule to assign $p(x)$
to a vertex $x$ when it is a normal child of its parent $u$ is somewhat more sophisticated: We have to consider two cases,
namely whether $u$ itself was a normal child or a pseudo leaf of its parent $u'$. These two cases are separately considered in
step 3.2.2 of the algorithm. We have provided some figures (figures \ref {first-step}   and \ref {second-step}) to help the
reader visualize these steps, to some extent. Later, in section
 \ref {algo-in-section}
we have also considered an example tree, and described in detail, how our
algorithm will assign position vectors to the vertices in the case of that
tree.

We use $K(u)$ to denote the
set of corners of $B(u)$. Also, if $v$ is a ``normal" child of its parent, the algorithm associates a number $J(v) \in \{1,2\}$ to
remember the axis along which shifting was done to get $p(v)$ (see step 3.2.2 of the algorithm and the comment about $J(v)$ in section
\ref {Comments}  for more details). We use
$\overline{\mathbf{0}}$ to denote the all zeroes vectors and $\overline{\mathbf{1}}$ to denote the all ones vector in $d$
dimensions.

~\\
\begin{tabbing}
\Large{\textbf{Algorithm}}~~~\=
~\\
\textbf{INPUT:}\ The rooted tree $T'$ obtained from $T$ in Section~\ref{Leaf-degree}, with a pseudo leaf chosen\\
 from $C(x)$ for each vertex $x$ with $L(x) =
 \emptyset$. (See section \ref {Basic Notation}).  \\
~\\
\textbf{OUTPUT:}\ Two functions \=$p: V(T) \rightarrow \mathbb R^d$ and $r: V(T) \rightarrow \mathbb R^+$ where $d = \lceil \log_2(\beta+2)\rceil$\\
~\\
\textbf{Step 1:} For the root $x$, $p(x)=\overline{\mathbf{0}}$ and $r(x)=1$.\\
\textbf{Step 2:} For the unique child $x'$ of $x$, $p(x')=\overline{\mathbf{1}}$ and $r(x')=\dfrac{r(x)}{8(|A(x)|+1)}=\dfrac{1}{8}$\\
        as $A(x) = \emptyset$, since in this case $C(x) = \{x'\}$ and $x'$ will be chosen as the pseudo\\ leaf, if it is not a leaf.\\
\textbf{Step 3:} Suppose $u$ is a \textbf{non-root} vertex for which $r(u)$ and $p(u)$ is\\
\hspace{5mm}already defined by the algorithm.\\
\hspace{5mm}\textbf{Step 3.1: (Defining $\mathbf{r(y)}$ for $\mathbf{y\in C(u)}$)} For each $y \in C(u)$ do:\\
\hspace{1cm}If $y$ is a leaf, then $r(y)=r(u)$\\
\hspace{1cm}Else $r(y)=\dfrac{r(u)}{8(|A(u)|+1)}$\\
~\\
\hspace{5mm}\textbf{Step 3.2: (Defining $\mathbf{p(y)}$ for $\mathbf{y\in C(u)}$)} For each $y \in C(u)$ do:\\
\hspace{1cm}Let $u'$ be the parent of $u$.\\
\hspace{1cm}Let $K'(u)=K(u) - B(u')$.\\
\hspace{15mm}\textbf{Step 3.2.1: Defining $p(y)$ for $y\in C(u) - A(u)$} {\bf (see figure \ref {first-step}) } \\
\hspace{2cm}Assign a point from $K'(u)$ to $p(y)$ such that no two vertices \\
\hspace{2cm}from $C(u)-A(u)$ are assigned the same position vector,\\
\hspace{2cm}i.e., for $y, y' \in C(u) - A(u)$ we have $p(y) \neq p(y')$ if $y \neq y'$.\\
\hspace{2cm}(\textit{See Lemma~\ref{Three in One} for the feasibility of this step}).\\
~\\
\hspace{15mm}\textbf{Step 3.2.2: Defining $p(y)$ for $y\in A(u)$}\\
\hspace{2cm} {\bf [First Case:] } If $u \notin A(u')$ \textbf{(i.e. if $u$ is a pseudo-leaf.) } \\
\hspace {2cm}   (See figure \ref {second-step} for the illustration of
  this step.)\\
\hspace{2.5cm}then let $q \in K'(u)-\{p(y): y \in C(u) - A(u)\}$.\\
\hspace{2.5cm}(\textit{See Lemma~\ref{Three in One} for the feasibility of this step}).\\
~\\
\hspace{2.5cm}Let $q^1$ be the other end-point of the edge $E_1(q)$.\\
\hspace{2.5cm}Let $A(u)=\{u_1,u_2,\ldots,u_t\}$ where $t=|A(u)|$.\\
\hspace{2.5cm}Let $p^i = q+(\frac{r(u)}{2})(1 +\frac{i}{t+1})\frac{q^1-q}{2r(u)} \ \ \forall \ i\in [t]$.\\
\hspace{2.5cm}($\frac{q^1-q}{2r(u)}$ is the unit vector along edge $E_1(q)$ from $q$ to $q^1$)\\
~\\
\hspace{2.5cm}\textbf{Shifting $\mathbf{p^i}$ along axis $2$ to get $\mathbf{p(u_i)}$}\\
\hspace{2.5cm}(Note that $d \geq 2$. Therefore we have at least 2 axes.)\\
\hspace{2.5cm}Define $p(u_i) = p^i(2,\dfrac{r(u)}{16(t+1)}) \ \ \forall \ i\in [t]$.\\
\hspace{2.5cm}Set $J(u_i) = 2\ \ \forall \ i\in [t]$.\\
\hspace{2.5cm}\\
\hspace{2cm}  {\bf [Second Case:] } Else if $u \in A(u')$ \textbf{(i.e. if $u$ is not a pseudo-leaf)}\\
\hspace{2.5cm}Let $q \in  K(u) \cap B(u')$.\\
\hspace{2.5cm}\textit{(See Lemma~\ref{Three in One} for the feasibility of this step)}\\
\hspace{2.5cm}Let $l = J(u)$\\
~\\
\hspace{2.5cm}Let $A(u)=\{u_1,u_2,\ldots,u_t\}$ where $t=|A(u)|$.\\
\hspace{2.5cm}Let $p^i = q + (\frac{r(u)}{2})(1 +\frac{i}{t+1})\frac{q^l-q}{2r(u)} \ \ \forall \ i\in [t]$.\\
\hspace{2.5cm}($\frac{q^l-q}{2r(u)}$ is the unit vector along edge $E_l(q)$ from $q$ to $q^l$)\\
\hspace{2.5cm}\\
\hspace{2.5cm}If $l=1$ then assign $j=2$ else assign $j = 1$ \\
\hspace{2.5cm}For each $u_i\in A(u)$, set $J(u_i) = j$.\\
\hspace{2.5cm}\\
\hspace{2.5cm}\textbf{Shifting $\mathbf{p^i}$ along $J(u_i) = j^{th}$ axis to get $\mathbf{p(u_i)}$}\\
\hspace{2.5cm}Define $p(u_i) = p^i(j,\dfrac{r(u)}{16(t+1)}) \ \ \forall \ i\in [t]$.
\end{tabbing}

\label{Comments}

\begin{figure}[!h]
\centering
\includegraphics[width=5in]{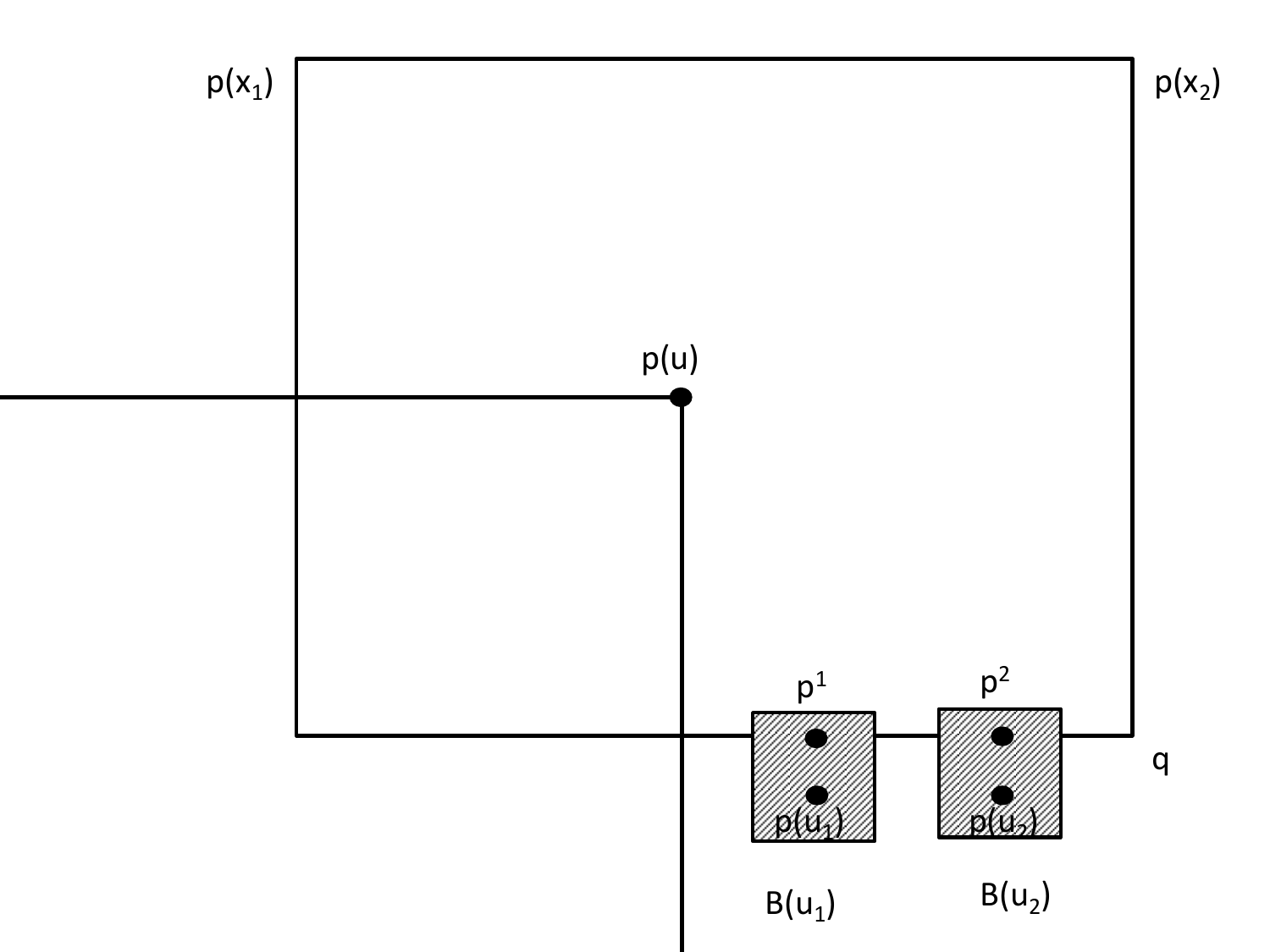}
\vspace{-5mm}
\caption
{
This figure is a continuation of Figure~\ref{second-step},  highlighting  $B(u)$ along with
 $B(u_i)$, $1\le i \le 2$. In this figure the shaded boxes represent $B(u_1)$ and $B(u_2)$.
Lemma~\ref{Two Corners} states that if a vertex is a normal child then exactly two
corner of its box lie inside the box of its parent. For example, in this figure each shaded box has
 exactly two corners inside the unshaded box corresponding to its parent.
\label{second-step-part-b}
}
\end{figure}

\subsection{Some Comments on the Algorithm}

Note that if a vertex $v$ is a normal child of its parent $u$ i.e. if $v\in A(u)$, then the algorithm in Step 3.2.2 assigns
the position $p(v)$ in one of the two possible ways depending on whether $u$ is a pseudo-leaf of its parent  or not. In both the cases, the
procedure is somewhat similar : We carefully select a corner $q$ of $B(u)$, then select a suitable edge $E_l(q)$ incident on
$q$, locate a suitable
position on $E_l(q)$ and then shift this position by distance $\frac{r(v)}{2} = \frac{r(u)}{16(|A(u)|+1)}$ along a chosen axis $j=J(v)\neq l$. (Note
that always $j,l \in \{1,2\}$.) \\
~\\
For the rest of the proof, for a vertex $v\in A(u)$, we say that $v$ is \textbf{attached} to the corner $q$ of $B(u)$ and the
edge $E_l(q)$ if the algorithm selects the corner $q$ and the edge $E_l(q)$ in order to find $p(v)$. (Note that according to the algorithm, all the
vertices of $A(u)$ get attached to the same corner and edge of $B(u)$.) \\
~\\
Let $A(u) = \{u_1,u_2,\ldots,u_t\}$ and let $v=u_i$. Let $v$ be attached to the corner $q$ of $B(u)$ and the edge $E_l(q)$.
Clearly $q=p(u)+ r(u).S_1$ for some $S_1\in {\cal S}$.
For presenting the proofs in later sections, it is important for us to be able to express $p(v)$ in terms of $q$.
To reach  the point $p(v)$ from the point $q$, we need to first move along the edge $E_l(q)$ and then
shift along the $j^{th}$ axis. Thus $p(v)$ and $q$ can be different only in $l^{th}$ and $j^{th}$ co-ordinates.
We can describe the co-ordinates of $p(v)$ in terms of the co-ordinates of $q$ precisely as follows:\\
If $k\notin \{l,j\}$, then $p(v)[k] = q[k]$.\\
For the $J(v)=j^{th}$ co-ordinate we have $p(v)[j] = q[j] + S_1[j].\frac{r(v)}{2}$.\\
For the $l^{th}$ co-ordinate we have $p(v)[l] = q[l] - S_1[l].\frac{r(u)}{2}. (1 + \frac{i}{t+1})$.\\

For example, consider the two dimensional case shown in Figure~\ref{fig-shifting}. Let the square $B(p,r)$ shown in this
figure correspond to $B(u)$ and let the child $v=u_i$ of $u$ be attached to the corner $q = p + r [1,-1]^T$ of $B(u)$.
That is, $S_1 = [1,-1]^T$ in this example.  Let
$E_1(q)$ be the edge to which $v$ is attached. That is $l = 1$ in  this example. Also $j = 2$.   Let $\epsilon =
\frac{r(u)}{2}. (1 + \frac{i}{t+1})$, where $t = |A(u)|$  and let $\delta = \frac {r(v)}{2}$.  Then clearly $z(2,\delta)$ in the figure,
corresponds to $p(v)$. The reader may want to verify that the first coordinate of $z(2, \delta)$ is given by $z(2,\delta)[1] =
z[1] = q[1] - S_1[1] \epsilon = q[1] - \epsilon =  p(u)[1] + r - \epsilon$. Similarly, $z(2,\delta)[2] = z[2] + S_1[2] \delta =
z[2] - \delta = q[2] - \delta = p(u) [2] - r - \delta$.

The ideas described in the above two paragraphs  is summarized in the following lemma :

\begin{lemma}
 Let $v\in A(u)$ and $|A(u)| = t$. Let $v$ be attached to the corner $q$ of $B(u)$ and the edge $E_l(q)$. Since $q\in K(u)$, we have
$q = p(u) + r(u).S_1$ for some $S_1\in {\cal S}$. If $\{e_1,e_2,\ldots,e_d\}$ is the canonical basis of $\mathbb{R}^d$, then
$p(v) = p(u) + r(u).S_1 + S_1[j].\frac{r(v)}{2}.e_j - S_1[l].\frac{r(u)}{2}.(1 + \frac{i}{t+1}).e_l$ for some $i\in [t]$ and
$j=J(v)\neq l$. As a consequence, note that $\rho(p(u),p(v)) = |p(u)[j]-p(v)[j]| = r(u) + \frac{r(v)}{2}$
\label{Co-ordinates}
\end{lemma}

\begin{lemma}
 Let $v\in A(u)$ be attached to corner $q\in K(u)$. Then no vertex in $C(u)-A(u)$ will be assigned to the corner $q$.
\label{No Other}
\end{lemma}
\begin{proof}
 If $u$ is a pseudo-leaf, then recalling Step 3.2.2 of the algorithm, we know that $q \in K'(u)-\{p(y)\ |\ y \in C(u) - A(u)\}$
 and therefore the lemma is true.\\
 If $u$ is not a pseudo-leaf, then recalling Step 3.2.2 of the algorithm, we know that $q \in  K(u) \cap B(u')$. So,
 $q\notin K'(u)$. But from Step 3.2.1 of the algorithm, each vertex in $C(u)-A(u)$ gets assigned to a corner from $K'(u)$.
 \hfill\qed
\end{proof}

\noindent {\bf Comment about $J(u)$:}  Note that the algorithm assigns a value from $\{1,2\}$ to all normal vertices of the
tree. (Thus the variable $l$ and $j$ in the algorithm always gets values from $\{1,2\}$.) If the parent $u'$ of the normal
vertex $u$ was a pseudo leaf, then $J(u) = 2$ (first case of step 3.2.2); else $J(u)$ is given the value $1$ or $2$ so that
$J(u) \ne J(u')$ (second case of step 3.2.2). The significance of $J(u)$ can be understood as follows: If $u$ is a normal
vertex, then $B(u)$ has exactly one edge inside $B(u')$ where $u'$ is the parent of $u$ (see Lemma~\ref{Two Corners} for a
formal proof of this statement).  This edge will be parallel to  axis 1 if $u'$ is a pseudo leaf. If $u'$ is also a normal
child of its own parent, then it will be parallel to the axis $l=J(u')$, which can be 1 or 2 as we have already seen. (Also see the comment after
Lemma \ref {Two Corners}.)  Thus
$J(u)$ is the smallest axis such that it is NOT parallel
to the edge of $B(u)$ that is inside $B(u')$.

\subsection{An Example of the Algorithm }
\label{algo-in-section}

\begin{figure}[!h]
\centering
\includegraphics[width=5.5in]{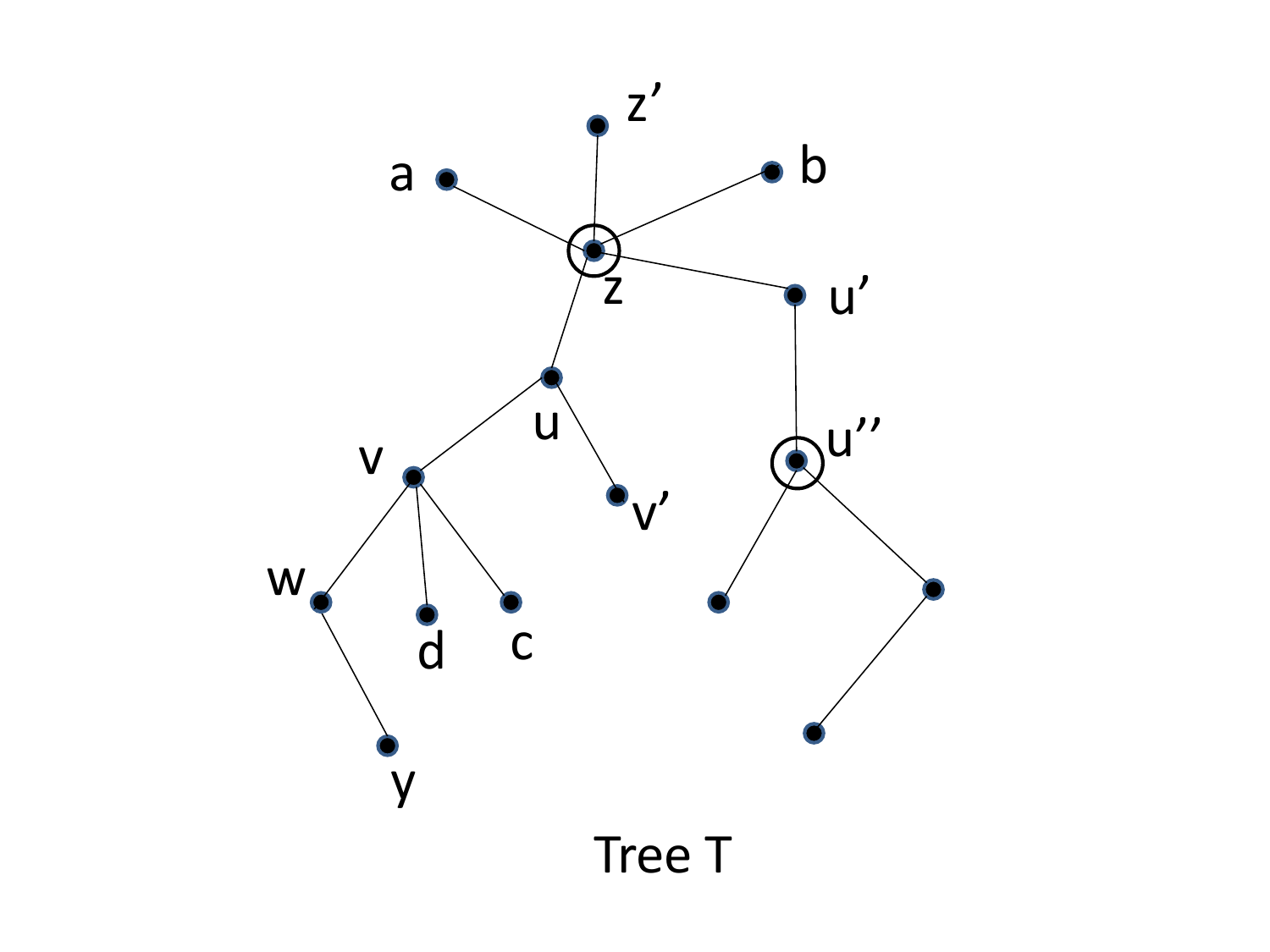}
\caption{In this figure we have a tree rooted at $z'$ as per Section~\ref{Leaf-degree}.
Here $\beta(T) = 2$. There are two pseudo leaves $z,u''$, which are marked
with the circles in the figure.
In Figure~\ref{algo-example-2} we show the boxes
for the vertices $z,u,v$ and $w$. For some other vertices just the position
is shown.
\label{algo-example-1}
}
\end{figure}

In Figure~\ref{algo-example-1} we have a tree $T$. Following the notation of Section~\ref{Leaf-degree} we have  $\alpha=3$.
Since there is only one vertex, namely $z$ with $3$ leaves adjacent to it, we select one of the leaves attached to it, say
$z'$ as the root. Clearly, $\beta = \alpha-1 = 2$ in this case. The algorithm of Section~\ref{Algorithm} computes position
vectors and radius for the vertices of this tree  in  $\lceil \log_2(\beta+2)\rceil = 2$ dimensions.

\begin{figure}[!h]
\centering
\includegraphics[width=5.5in]{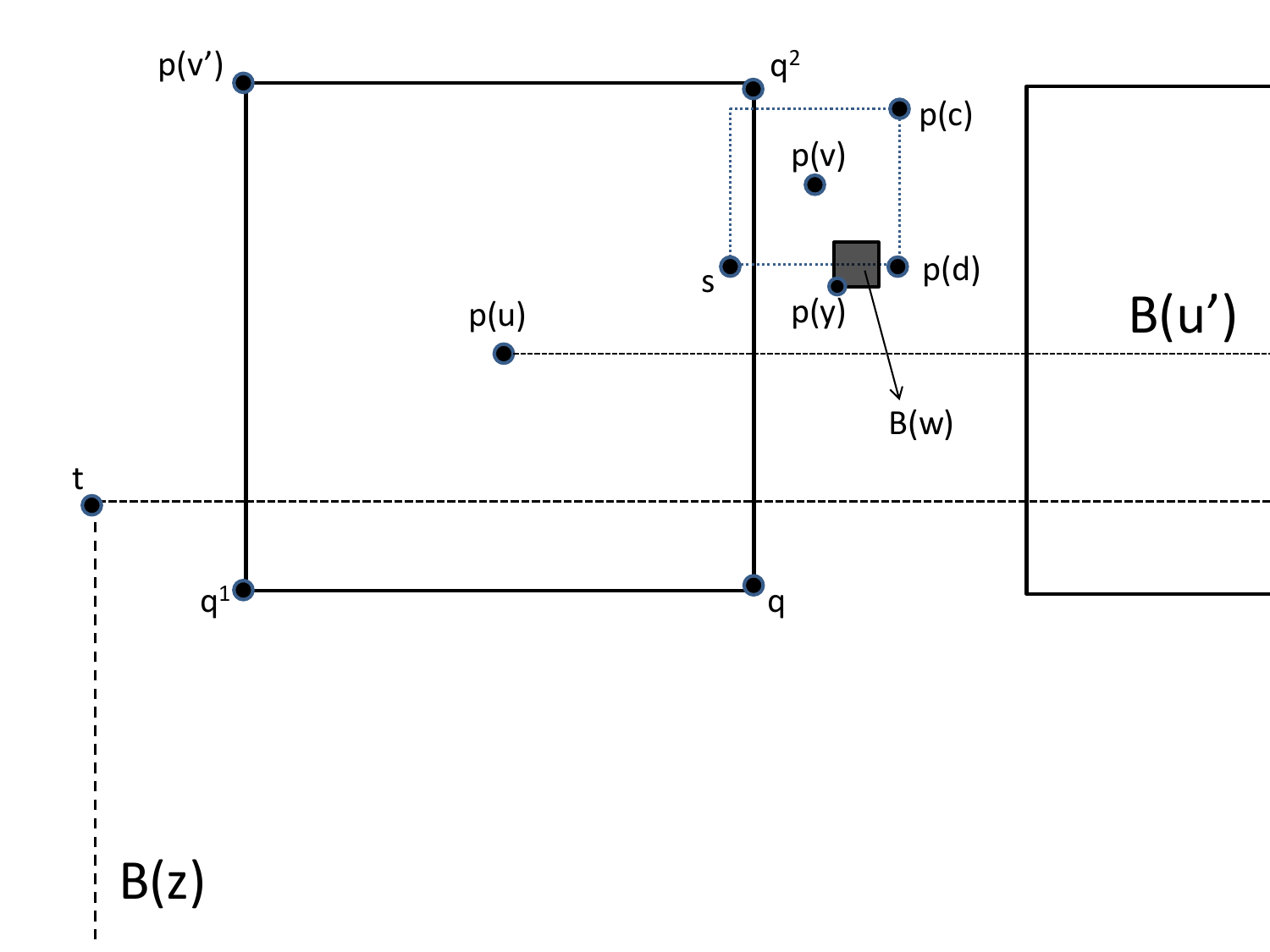}
\caption{In this figure we show the boxes of $z,u,v,u'$ and $w$.
Also the positions of $v',c,d,y$ are also marked. Assuming that X axis is
the first axis  and Y axis is the second axis, then the diagram indicates that $J(u)=2$.
Hence the algorithm will shift along the 1st axis to place $v$.
\label{algo-example-2}
}
\end{figure}

We will describe here how the algorithm of Section~\ref{Algorithm} computes the position vectors for some of the vertices of
this tree. In Figure~\ref{algo-example-2}, we have illustrated (approximately) how these positions can be marked on the
$2$-dimensional plane. For some vertices, we have shown the corresponding balls (which are squares in $2$-dimensions).
Unfortunately, it is difficult to draw the squares for  all the vertices of the tree: That is why we have selected only few
vertices to illustrate.  Moreover, the figure is not drawn to scale since the actual radius of the squares, decreases quite
fast, as the distance of the corresponding vertices from the root increases. Thus these illustration should be considered only
as an aid to the reader to intuitively visualize the procedure.

Algorithm first assigns the position $(0,0)$ and a radius of $1$ to the root $z'$. In other words, the algorithm assigns a
unit square to $z'$ centered at $(0,0)$. By step 2 of the algorithm, $z$, which is the only child of $z'$, will be assigned
the position corresponding to the top right hand corner of the square of $z'$. In Figure~\ref{algo-example-2} we have not
shown $B(z')$. The square for $z$,  $B(z)$ is only partially shown and is drawn  with dotted lines at the bottommost. Now, $z$
has $4$  children: $u,u',a,b$.
  By step 3.2.1, the $2$ leaves $a$ and $b$  will be assigned
positions which are corners of $B(z)$ and outside $B(z')$.
We consider now how the algorithm will assign positions to $u,u'$.
 As $u$ and $u'$ are  normal children  of $z$,
we are at step 3.2.2 of the algorithm. Moreover, since
the parent of $u$ namely $z$ is a pseudo leaf of its parent $z'$, we are in the first case of step 3.2.2.  Let $t$ be the
corner of $B(z)$ that the algorithm selected to attach
 $u$ and $u'$  to. Then $E_1(t)$ is the edge of $B(z)$ that $u$ and
$u'$  will get attached  by first case of
step 3.2.2. The algorithm computes the position for $u$ and $u'$
 by first moving
a carefully calculated distance along $E_1(t)$ and then
 shifting along the 2nd axis: the position of $u$, $p(u)$ and its
box is shown in Figure~\ref{algo-example-2}.
In a similar way, $u'$ will be assigned a position at the same `height' as
$p(u)$, but a little more to the right of $p(u)$. The box $B(u')$ is shown
partially in figure. (Caution: The figure is not to scale, and the reader
should not worry that the distances are not matching with what the algorithm
prescribes.)
 Note that $J(u)=J(u')=2$.

Now we have  to place the children of $u$ : $v$ and $v'$.
Since $v'$ is a leaf, it will be placed first, by step 3.2.1,
 at one of the external corners
of $B(u)$, see the figure.
Now  $v$ is placed by the  second case of Step 3.2.2, as $v$ is a normal child of $u$ and $u$ is not a pseudo leaf of $z$. By
Lemma~\ref{Two Corners}, $B(u)$ has exactly 2 corners inside $B(z)$: This
can be seen clearly from the figure. Let $q$ be one of those corners.
We then traverse along
the edge $q-q^{l}$ where $l=J(u)=2$ and $q^l$ is the corner of $B(u)$ which differs from $q$ only in the $l^{th}$
coordinate. We traverse along the edge $q-q^{2}$ a distance
more than $\frac{r(u)}{2}$ and then shift along an axis different
from  $J(u)$ to get the position for $v$. As $J(u)=2$ we will be shifting
 along the first axis in this case to get the position
for $v$. The position $p(v)$ and $B(v)$ is shown in figure. Also note
that $J(v) = 1$, since $l = J(u) = 2$.

Two of the children $c,d$ of $v$ are leaves and they are placed as shown in
figure at the external corners of $B(v)$. Since $w$ is a normal child of
$v$ and $v$ itself is a normal child of its parent $u$, the second case of
step 3.2.2 will be assigning the position for $w$. Let $s$ be the corner of
$B(v)$ to which $w$ gets attached to: Since $J(v)=1$,  $w$ will be attached to
the edge $E_1(s)$.
Then the box $B(w)$ corresponding to
$w$ will be placed approximately as shown in figure.

Finally, the vertex $y$ which is a leaf child of $w$ will be assigned to an
external corner of $B(w)$ as shown in figure.

\subsection{Correctness of the Algorithm}
\label{Correctness}

In this section we provide the lemmas to establish that the algorithm given in Section~\ref{Algorithm} does not get stuck at
any step and thus assigns to each vertex a radius and a position vector.

After reading the statements of Lemma \ref {One Corner} and  \ref {Two Corners}, the reader is requested to
verify this in the 2-dimensional cases given in figure \ref {first-step} (for the case of
Lemma \ref{One Corner})  and figure \ref {second-step} followed by figure \ref {second-step-part-b} (for the case of
Lemma \ref {Two Corners}).  In the
2-dimensional case, these Lemmas are intuitive and immediate, but for higher dimensional case,  we need a rigorous proof
which is provided below.

\begin{lemma}
Let $u'$ be the parent of $u$. If $u \in C(u')-A(u')$ i.e. $u$ is either a leaf or a pseudo-leaf , then exactly 1 corner of
$B(u)$ lies inside $B(u')$. Moreover, let $p(u)$ be the corner of $B(u')$ given by $p(u')+r(u').S_1$ where $S_1\in {\cal S}$.
Then the corner of $B(u)$ which is inside $B(u')$ is given by $p(u)-S_1.r(u)$.
\label{One Corner}
\end{lemma}
\begin{proof}
Let ${\cal S} = \{-1,+1\}^d$. We shift co-ordinates so that $p(u') = \overline{\textbf{0}}$. Since $u \in C(u')-A(u')$, by
Step 3.2.1 of the algorithm, $p(u)$ gets assigned to one of the corners of $B(u')$. So, we have $p(u) = p(u')+r(u')S_1 =
r(u')S_1$ for some $S_1\in {\cal S}$. Now consider  a general corner $q$ of $B(u)$. Since $q\in K(u)$, we know that $q = p(u)
+ r(u)S_2$ for some $S_2\in {\cal S}$. So, $q = r(u')S_1 + r(u)S_2$. For $q$ to be inside $B(u')$, we require
 $r(u') > \rho(p(u'),q)$ i.e. $r(u') > \rho(\overline{\textbf{0}},r(u').S_1 + r(u).S_2)$. We know that $r(u) > 0$. If there is any
 co-ordinate position $n\in [d]$ such that $S_1[n] = S_2[n]$, then $\rho(p(u'),q) = r(u') + r(u) > r(u')$.  So, we infer that
for $q$ to be inside $B(u')$, it is necessary that $S_2 =-S_1$. Also, if $S_2 =-S_1$, then $\rho(p(u'),q) = |r(u') - r(u)| =
r(u') - r(u) < r(u')$. Therefore, $q$ is inside $B(u')$ if and only if $S_2 =-S_1$. Clearly, given $S_1\in {\cal S}$ there is
a unique $S_2\in {\cal S}$ such that $S_2 =-S_1$. So, there is exactly one corner of $B(u)$ inside $B(u')$ and it is given by
$p(u)-S_1.r(u)$.
\hfill\qed
\end{proof}

\begin{lemma}
Let $u'$ be the parent of $u$. If $u \in A(u')$ i.e. if $u$ is a normal child of $u'$, then exactly $2$ corners of $B(u)$ lies
inside $B(u')$. Further, these $2$ corners form an edge of $B(u)$. Let $u$ be attached to corner $q$ of $B(u')$ and edge
$E_l(q)$ where $q = p(u') + r(u').S_1$ with $S_1\in {\cal S}$. Then the $2$ corners of $B(u)$ inside $B(u')$ are given by
$p(u) + r(u).S'$ and $p(u) + r(u).S''$ with $S',S''\in {\cal S}$ where $S'=-S_1$  and $S''$ is the string which differs from
$S_1$ in all co-ordinates except the $l^{th}$ co-ordinate.
\label{Two Corners}
\end{lemma}
\begin{proof}
Let ${\cal S} = \{-1,+1\}^d$ and $\{e_1,e_2,\ldots,e_d\}$ be the canonical basis of $\mathbb{R}^d$. Let $|A(u')|=t$. Without
loss of generality, let $p(u') = \overline{ \textbf{0}}$. Since $u\in A(u')$, let $u$ be attached to the corner $q\in K(u')$
and edge $E_l(q)$. Let $q = p(u')+ r(u').S_1 = r(u').S_1$ where $S_1\in {\cal S}$. Thus by Lemma~\ref{Co-ordinates}, $p(u) = q
+ S_1[j].\frac{r(u)}{2}.e_j - S_1[l].\frac{r(u')}{2}. (1 + \frac{i}{t+1}).e_l = r(u').S_1 + S_1[j].\frac{r(u)}{2}.e_j -
S_1[l].\frac{r(u')}{2}. (1 + \frac{i}{t+1}).e_l$ \ for some $i\in [t]$ and $j\neq l$. Consider any general corner $s$ of
$B(u)$. Then $s = p(u) + r(u). S_2$ for some $S_2\in {\cal S}$. So,
$$s = r(u').S_1 + S_1[j].\frac{r(u)}{2}.e_j - S_1[l].\frac{r(u')}{2}. (1 + \frac{i}{t+1}).e_l + r(u). S_2$$

Note that   $s$ is  inside $B(u')$,
if and only if  $r(u') > \rho(p(u'),s)$ i.e. $r(u') > \rho(\overline{\textbf{0}}\ ,s)$:  in other words, $|s[n]| < r(u'), 1 \le n \le d$.

\begin{enumerate}
\item $\mathbf{|s[l]| < r(u')}$\\
$|s[l]| = |S_1[l].r(u') - S_1[l].(\frac{r(u')}{2})(1 +\frac{i}{t+1}) + S_2[l]r(u)| \ \leq |r(u') - (\frac{r(u')}{2})(1 +\frac{i}{t+1})| +
|S_2[l]r(u)| = |(\frac{r(u')}{2})(1 -\frac{i}{t+1})| + |r(u)|$. Now, $i\in [t]$ and $r(u) = \frac{r(u')}{8(t+1)} \leq \frac{r(u')}{8}$.
Therefore, $|s[l]| \leq  |(\frac{r(u')}{2})(1 -\frac{i}{t+1})| + |r(u)| \leq |\frac{r(u')}{2}| + |\frac{r(u')}{8}| < r(u')$.
So, independent of value of $S_2[l]$, we get that $|s[l]| < r(u')$ .\\
\item $\mathbf{|s[j]|<r(u')}$ \textbf{if and only if} $\mathbf{S_2[j]=-S_1[j]}$\\
If $S_1[j] = S_2[j]$, then $|s[j]| \geq r(u') + \frac{3r(u)}{2} > r(u')$. On the other hand, if  $S_1[j] \neq S_2[j]$ and $s[j] = r(u') -
\frac{r(u)}{2} < r(u')$.\\
\item $\mathbf{|s[k]|<r(u')}$ \textbf{if and only if} $\mathbf{S_2[k]=-S_1[k] \ \forall \ k\in [d], \ k\notin\{l,j\}}$\\
Since $r(u) > 0$, if  $S_1[k]\neq S_2[k]$ we get  $|s[k]| = r(u') - r(u) < r(u') \ \ \forall\  k \notin \{l,j\}$. On
 the other hand, if $S_1[k] = S_2[k]$, we get $|s[k]| = r(u') + r(u) > r(u')  \ \ \forall\  k \notin \{l,j\}$.  \\
\end{enumerate}
From the above we can conclude that for given $S_1\in {\cal S}$, the corner $s = p(u) + S_2.r(u)$ is inside $B(u')$ if and
only if $S_2[n]\neq S_1[n] \ \ \forall n\neq l$. We infer that the strings $S',S''$  ( where $S'=-S_1$  and $S''$ is the
string which differs from $S_1$ in all co-ordinates except the $l^{th}$ co-ordinate ) correspond to the two corners of $B(u)$
inside $B(u')$.\hfill\qed
\end{proof}

~~~~~~

\noindent {\bf Comment:} Note that the edge between the two corners of $B(u)$ that are inside $B(u')$ is parallel to the
$l^{th}$ axis.  This is because these two corners, namely $p(u) + r(u).S'$ and $p(u) + r(u).S''$ differs in only the $l^{th}$
co-ordinate, by the statement of Lemma \ref {Two Corners}. If $u'$ also was a normal child, then $l=J(u')$ (which can be $1$
or $2$)  by the second case of step 3.2.2 of the  algorithm. If $u'$ was a pseudo leaf of its parent then $l=1$, by the first
case of step 3.2.2 of the algorithm.

~~~~

\begin{lemma}
 Let $u$ be any non-root non-leaf vertex. Let $K'(u) = K(u) - B(u')$ where $u'$ is parent of $u$. Then,
 if $u$ is a pseudo-leaf  we have $|K'(u)|>\beta \geq|C(u) - A(u)|$ and if $u$ is a normal child we have $|K'(u)|\geq \beta\geq
|C(u) - A(u)|$
\label{More Outside}
\end{lemma}
\begin{proof}
Note that $|C(u) - A(u)|$ is the number of children of $u$ that are not ``normal''. Recall from Section~\ref{Leaf-degree}
that $\beta \geq $ 1. If $u$ has a child which is a leaf, then $|C(u) - A(u)|\leq \beta$. Else, $|C(u) - A(u)| = 1 \leq \beta$.\\
 If $u$ is a pseudo-leaf of its parent $u'$, then by Lemma~\ref{One Corner} and recalling that
 $d = \lceil \log_2(\beta+2)\rceil$, $|K'(u)| = |K(u)| - 1=2^d - 1 \geq \beta + 1 > \beta$.\\
If $u$ is not a pseudo-leaf of its parent $u'$, then by Lemma~\ref{Two Corners} and recalling that $d = \lceil
\log_2(\beta+2)\rceil$, $|K'(u)| = |K(u)| - 2= 2^d - 2 \geq \beta$.\hfill\qed
\end{proof}

\begin{lemma}
 The algorithm given in the previous section runs correctly i.e. it does not get stuck at any of the three points
 which give a reference to this lemma. Therefore each vertex is assigned a position vector and a radius when the algorithm terminates.
\label{Three in One}
\end{lemma}
\begin{proof}
By Lemma~\ref{More Outside}, we know that $|K'(u)|\geq |C(u) - A(u)|$. So, distinct points in $C(u)-A(u)$ can be assigned
distinct points from $K'(u)$ in Step 3.2.1 of the algorithm\\
If $u$ is a ``pseudo-leaf'' of its parent $u'$, then by Lemma~\ref{More Outside}, $|K'(u)| > \beta \geq |C(u) - A(u)|$. So,
$( K'(u)-\{p(y): y \in C(u) - A(u)\} ) \neq \emptyset$ as is required in the ``if'' part if Step 3.2.2 of the algorithm.\\
If $u$ is not a ``pseudo-leaf'' of its parent $u'$, then by Lemma~\ref{Two Corners}, there are two corners of $B(u)$ inside
$B(u')$ and therefore $K(u) \cap B(u') \neq \emptyset$ as is required in ``else'' part of Step 3.2.2 of the
algorithm.
\hfill\qed
\end{proof}

\subsection{$T$ is the Intersection Graph of the family $\{ B(u)\ |\ u\in V(T) \}$}
\label{First Step}

Recall that for every vertex $u$ in $V(T)$, $B(u) = B(p(u),r(u))$ where $p(u)$ and $r(u)$ are the position vector and radius
computed by the algorithm. In this section we show that $T$ is the intersection graph of the family $\{ B(u)\ |\ u\in V(T) \}$.

\begin{lemma}
Let $uv\in E(T)$. Then, $B(u) \cap B(v)\neq \emptyset$.
\label{Edge}
\end{lemma}
\begin{proof}
Immediate from Lemma~\ref{One Corner} and Lemma~\ref{Two Corners}: If $u$ is the parent of $v$, at
least one corner of $B(v)$ is inside $B(u)$.
\hfill\qed
\end{proof}

\begin{lemma}
 Let $v,w$ be children of $u$. Then $S(v) \cap S(w) = \emptyset$
\label{Superball of children}
\end{lemma}
\begin{proof}
We have the following 3 cases: \\
\textbf{Case 1} : \\
Both $v,w\in C(u)-A(u)$. Then both $v$ and $w$ must be leaves. Thus, $R(v) = R(w) = r(u) = r(w) = r(v)$. By Step 3.2.1 of the
algorithm,
$p(v)$ and $p(w)$ will ge placed at distinct corners of $B(u)$. So $\rho(p(v),p(w)) = 2r(u) = R(v) + R(w)$ and hence recalling
  that $S(v)$ and $S(w)$ are open balls,  $S(v) \cap S(w) = \emptyset$.\\
\textbf{Case 2 }: \\
Both $v,w\in A(u)$ and let $|A(u)| = t$. So, $r(v) = \frac{r(u)}{8(t+1)} = r(w)$ and $R(v)=2r(v)$ and $R(w)=2r(w)$. Following
terminology of Section~\ref{Comments}, let $q\in K(u)$ be the corner and let $E_l(q)$ be the edge to which $v$ and $w$ are
attached . Let $q= p(u) + S_1.r(u)$ where $S_1\in {\cal S}$. Applying Lemma~\ref{Co-ordinates}, we  have $p(v) = q +
S_1[j].\frac{r(v)}{2}.e_j - S_1[l].\frac{r(u)}{2}. (1 + \frac{i}{t+1}).e_l$ for some $i\in [t]$ and $p(w) = q +
S_1[j].\frac{r(w)}{2}.e_j - S_1[l].\frac{r(u)}{2}. (1 + \frac{k}{t+1}).e_l$ for some $k\in [t]$. Note that $j\neq l$. Also,
$v$ and $w$ get assigned to distinct points and hence $|i-k|\geq 1$. Recalling that $r(v) = r(w)$, we see that $p(v)$ and
$p(w)$ differ only in $l^{th}$ co-ordinate.
So, $\rho(p(v),p(w)) = |p(v)[l] - p(w)[l]| = |(i-k).\frac{r(u)}{2(t+1)}| \geq \frac{r(u)}{2(t+1)} = R(v) + R(w)$.
 Hence recalling that $S(v)$ and $S(w)$ are open balls,  $S(v) \cap S(w) = \emptyset$.\\
\textbf{Case 3} : \\
Let $v\in C(u)-A(u)$ and $w\in A(u)$. Let $|A(u)|=t$. So, $r(w) = \frac{r(u)}{8(t+1)} < \frac{r(u)}{8}$ and $R(w) = 2r(w) <
r(u)$. If $v$ is a leaf, then $r(v)=R(v)=r(u)$. If $v$ is a pseudo-leaf, then $R(v) = 2r(v)$ and $r(v) = \frac{r(u)}{8(t+1)}
\leq \frac{r(u)}{8}$. In either case, $R(v) \leq r(u)$. We translate the co-ordinates so that $p(u) = \overline{\textbf{0}}$.
By Step 3.2.1 of the algorithm, $v$ gets assigned to a corner $q_1\in K(u)$. So, $p(v) = p(u)+S_2.r(u)=S_2.r(u)$ where $S_2\in
{\cal S}$. Let $q\in K(u)$ be the corner and $E_l(q)$ be the edge to which $w$ is attached. So, $q= p(u)+S_1.r(u) = S_1.r(u)$, for
some $S_1 \in  {\cal S}$.
From Lemma~\ref{Co-ordinates}, we have
$p(w) = S_1.r(u) + S_1[j].\frac{r(w)}{2}.e_j - S_1[l].\frac{r(u)}{2}. (1 + \frac{i}{t+1}).e_l$ for some $i\in [t]$ and $j\neq l$. \\
~\\
\textit{Claim.} $\mathbf{q\neq q_1}$\\
Recall that at Step 3.2.1 of the algorithm the corner $q_1$ of $B(u)$ given to $p(v)$ is from $K'(u)$ since $v\in C(u)-A(u)$.
Now, we know that $w\in A(u)$. In Step 3.2.2 of the algorithm, if $u$ is a pseudo-leaf, then $q$ is chosen from $K'(u)-\{p(y):
y \in C(u) - A(u)\}$. If $u$ is not a pseudo-leaf, then $q$ is chosen from $K(u) \cap B(u')$ which is disjoint from $K'(u)$.
So, in both the cases we get that $q\neq q_1$.
$\bbox$

~\\
 In view of the above claim, $q$ and $q_1$ (and therefore $S_1$ and $S_2$) differ in at least one co-ordinate say $k$. If $k\notin \{l,j\}$,
 then considering distance along $k^{th}$ co-ordinate, $\rho(p(v),p(w)) \geq |p(v)[k]-p(w)[k]| = 2r(u) \geq R(v) + R(w)$. If $k=j$, then
 considering distance along $j^{th}$ co-ordinate , we have $\rho(p(v),p(w)) \geq |p(v)[j]-p(w)[j]| = 2r(u) + \frac{r(w)}{2} > 2r(u) \geq R(v) + R(w)$.
 If $k=l$, the considering distance along $l^{th}$ co-ordinate we get
 $\rho(p(v),p(w))\geq |p(v)[l]-p(w)[l]| = 2r(u) - \frac{r(u)}{2}(1 + \frac{i}{t+1}) = r(u) + (r(u) - \frac{r(u)}{2}(1 + \frac{i}{t+1}))
 = r(u) + (\frac{r(u)}{2}(1 - \frac{i}{t+1})) \geq r(u)
 + \frac{r(u)}{2(t+1)} $ as $i\in [t]$. Hence, recalling that $r(w)=\frac{r(u)}{8(t+1)}$ we get $\rho(p(v),p(w)) \geq r(u) + \frac{r(u)}{2(t+1)}
 > r(u) + 2r(w) = R(v) + R(w)$. Hence $S(v) \cap S(w) = \emptyset$.
 \hfill\qed
\end{proof}

\begin{lemma}
If $v$ is a child of $u$, then $S(v) \subset S(u)$
\label{Child in Parent}
\end{lemma}
\begin{proof}
 We have the following 2 cases depending on what type of child $v$ is.\\
\textbf{Case 1 }:\\
$v\in C(u)-A(u)$. Recall that if $v$ is a leaf then $R(v)=r(v)=r(u)$ and if $v$ is a pseudo-leaf then $R(v)=2r(v)<r(u)$. In
both cases, we have $r(u)\geq R(v)$ Also note that $R(u)=2r(u)$ as $u$ is not a leaf. By Step 3.2.1 of the algorithm
$\rho(p(u),p(v)) = r(u)$ . Let $p(y)$ be any point in $S(v)$. Then, $\rho(p(v),p(y)) < R(v)$. By triangle inequality,
$\rho(p(u),p(y)) \leq \rho(p(u),p(v)) +
\rho(p(v),p(y)) < r(u) + R(v) \leq 2r(u) = R(u)$ . Hence $S(v) \subset S(u)$.\\
\textbf{Case 2 }:\\
$v\in A(u)$. Let $|A(u)|=t$. Then, $r(v) = \frac{r(u)}{8(t+1)} < \frac{r(u)}{8}$ and thus $R(v) = 2r(v) < \frac{r(u)}{4}$.
Also note that $R(u)=2r(u)$ as $u$ is not a leaf. By Lemma~\ref{Co-ordinates}, $\rho(p(u),p(v)) = r(u) + \frac{r(v)}{2}$ . Let
$p(y)$ be any point in $S(v)$. Then, $\rho(p(v),p(y)) < R(v)$. By triangle inequality, $\rho(p(u),p(y)) \leq \rho(p(u),p(v)) +
\rho(p(v),p(y)) < r(u) + \frac{r(v)}{2} + R(v) < r(u) + \frac{r(u)}{16} + \frac{r(u)}{4} < 2r(u) = R(u)$. Hence $S(v) \subset
S(u)$.
\hfill\qed
\end{proof}

\begin{lemma}
 Let $u$ be any vertex of $T'$. For every descendant $v$ of $u$, $B(v) \subseteq S(u)$.
\label{Every descendant}
\end{lemma}
\begin{proof}
Let the depth of a vertex $x$ in $T'$ be the number of vertices in the path from the root to $x$. We prove the lemma by
induction on the depth. Let $D$ be the maximum depth among all vertices of $T'$. Clearly any vertex of depth $D$ must be a
leaf. Now, if $u$ is a leaf then its only descendant is itself and the lemma holds trivially. So, we infer that the lemma
holds for all vertices of depth $D$.We take this as the base case of the induction. Now suppose that the lemma holds true for
all vertices of depth greater than $h$ where $1 \le h<D$. Now, let $u$ be a vertex of depth $h$. If $u$ is a leaf then the lemma
holds trivially. Otherwise, let the children of $u$ be $\{u_1,u_2,\ldots,u_m\}$. Clearly, depth of $u_i$ equals $(h+1)$
$\forall \ i\in [m]$. Now, let $v$ be a descendant of $u$ such that $v\neq u$. If $v=u_j$ for some $j\in [m]$, then by
Lemma~\ref{Child in Parent} we have $B(v) \subseteq S(v) \subset S(u)$. Else, $v\neq u_j \ \forall \ j\in [m]$. Then $v$ is a
descendant of $u_k$ for some $k\in [m]$. Recalling that depth of $u_k$ equals $(h+1)$ and using induction hypothesis, we have
$B(v) \subseteq S(u_k)$.
 By Lemma~\ref{Child in Parent}, $S(u_k) \subset S(u)$. Therefore, $B(v) \subseteq S(u_k) \subset S(u)$.
 \hfill\qed
\end{proof}

\begin{lemma}
 Let $u$ be the parent of $v$ and $v$ be the parent of $w$. Then $B(u) \cap S(w) = \emptyset$
\label{Grandchild}
\end{lemma}
\begin{proof}
Let ${\cal S}=\{-1,+1\}^d$. We shift co-ordinates so that $p(u) = \overline {\textbf{0}}$. Note that $v$ cannot be a leaf. So, we have
the following four cases.\\
\textbf{Case 1 : $\mathbf{v}$ is a pseudo-leaf}\\

\textbf{Case 1.1 : $\mathbf{w\in C(v)-A(v)}$ :}\\
So $p(v)$ is a corner of $B(u)$ and is hence given by $p(v) = p(u) + r(u).S_1 = r(u).S_1$ for some $S_1\in {\cal S}$. By Step
3.2.1 of the algorithm, we know that $p(w)$ belongs to $K'(v)$. Thus $p(w) = p(v) + r(v).S_2 = r(u).S_1 + r(v).S_2$ for some
$S_2\in {\cal S}$. Since $p(w) \in K'(u)$ we can infer that $S_2 \neq -S_1$ by Lemma~\ref{One Corner}. That is, there is some
$k\in [d]$ such that $S_1[k] = S_2[k]$. Considering distance along $k^{th}$ co-ordinate, we have $\rho(p(u),p(w)) =
\rho(\overline {\textbf{0}} , r(u).S_1 + r(v).S_2) = r(u) + r(v)$. Now, if $w$ is a leaf then $R(w) = r(w) = r(v)$. Else, if
$w$ is a pseudo-leaf, then $r(w)\leq \frac{r(v)}{8}$ and $R(w) = 2r(w)$.
In both cases, we are assured that $r(v)\geq R(w)$.Hence, $\rho(p(u),p(w)) =  r(u) + r(v)\geq r(u) + R(w)$. Therefore $B(u)\cap S(w)=\emptyset$.\\

\textbf{Case 1.2 : $\mathbf{w\in A(v)}$ :}\\
Let $|A(v)|=t$. Note that from Step 3.1, $r(w) = \frac{r(v)}{8(t+1)} < \frac{r(v)}{8}$ and thus $r(v)> 2r(w)=R(w)$. Also as in
the previous case, $p(v) = p(u) + r(u).S_1 = r(u).S_1$ for some $S_1\in {\cal S}$. Suppose that $w$ gets attached to corner
$q\in K(v)$ and edge $E_l(q)$. (See Section~\ref{Comments}). Clearly $q = p(v) + S_2.r(v) = r(u).S_1 + r(v).S_2$ for some
$S_2\in {\cal S}$. By Lemma~\ref{Co-ordinates}, we have $p(w) = r(u).S_1 + r(v).S_2 - \frac{r(v)}{2}(1+
\frac{i}{t+1}).S_2[l].e_l + S_2[j].\frac{r(w)}{2}.e_j$ for some $i\in [t]$ and $j\neq l$. By Step 3.2.2 of the algorithm, we
know that $q \in K'(v)-\{p(y)\ |\ y \in C(v) - A(v)\}$. Also from Lemma~\ref{One Corner} we know that the only corner of
$B(v)$ that is inside $B(u)$ is given by $p(v) - S_1.r(v)$ and since $q \in K'(v)$ is outside $B(u)$ we can infer that $S_2
\neq -S_1$.  i.e. there is $k\in [d]$ such that $S_1[k]= S_2[k]$. If $k\notin \{l,j\}$, then considering the distance along
$k^{th}$ co-ordinate $\rho(p(u),p(w))=\rho(\mathbf{\overline {0}}, p(w)) \geq r(u)+ r(v) >r(u) + R(w)$. If $k=j$, then
$\rho(p(u),p(w)) \geq r(u)+ r(v) + \frac{r(w)}{2} > r(u) + r(v) > r(u) + R(w)$. If $k=l$, then $\rho(p(u),p(w)) \geq r(u) +
r(v) -\frac{r(v)}{2}(1 + \frac{i}{t+1}) = r(u)+ \frac{r(v)}{2}(1- \frac{i}{t+1})$. Note that $i\in [t]$ implies that
$\frac{r(v)}{2}(1- \frac{i}{t+1}) \ge  \frac{r(v)}{2(t+1)} > 2r(w) = R(w)$.
So, $\rho(p(u),p(w)) \geq r(u)+ \frac{r(v)}{2}(1- \frac{i}{t+1}) \geq r(u) + R(w)$. Therefore $B(u)\cap S(w)=\emptyset$.\\

~\\

\textbf{Case 2 : $\mathbf{v}$ is not  a pseudo-leaf i.e. $\mathbf{v\in A(u)}$}\\

\textbf{Case 2.1 : $\mathbf{w\in C(v)-A(v)}$ :}\\
If $w$ is leaf, then $R(w)=r(w)=r(v)$. If $w$ is a pseudo-leaf, then by Step 3.1, we have $r(w) \leq \frac{r(v)}{8}$ and
therefore $R(w)=2r(w)<r(v)$. In both cases, $R(w)\leq r(v)$. Let $|A(u)|=t$. Suppose that $v$ gets attached to corner $q\in
K(u)$ and edge $E_l(q)$. (See Section~\ref{Comments}). Clearly $q=p(u)+S_1.r(u)=S_1.r(u)$ for some $S_1\in {\cal S}$. By
Lemma~\ref{Co-ordinates}, we have $p(v) = r(u).S_1 - \frac{r(u)}{2}(1+ \frac{i}{t+1}).S_1[l].e_l + S_1[j].\frac{r(v)}{2}.e_j$
for some $i\in [t]$ and a suitably selected $j\neq l$. Also, Step 3.2.1 of the algorithm assigns $p(w)$ to a corner of $B(v)$
and hence $p(w)=p(v)+S_2.r(v)$ for some $S_2\in {\cal S}$. That is, $p(w) = r(v).S_2 + r(u).S_1 - \frac{r(u)}{2}(1+
\frac{i}{t+1}).S_1[l].e_l + S_1[j].\frac{r(v)}{2}.e_j$. By Step 3.2.1 of the algorithm, we know that $p(w)\in K'(v)$. Also, by
Lemma~\ref{Two Corners}, we know that the 2 corners of $B(v)$ that are in $B(u)$ are given by $p(v) + r(v).S'$ and $p(v) +
r(v).S''$  where the strings  $S'$ and $S''$ are such that $S' = -S_1$ and $S''$ is the string which differs from $S_1$ in all
co-ordinates other than $l^{th}$ co-ordinate. It follows that $\exists \ m\in [d], m\neq l$ such that $S_1[m]=S_2[m]$ since
$p(w) = p(v) + S_2.r(v)$ is a corner of $B(v)$ outside $B(u)$. If $m=j$, then $\rho(p(u),p(w)) = \rho(\mathbf{\overline {0}},
p(w)) \geq r(v) + r(u) + \frac{r(v)}{2} > r(v)+r(u) \geq R(w)+r(u)$. If $m\neq j$, then $\rho(p(u),p(w)) \geq r(v)+r(u) \geq R(w)+r(u)$.
Therefore $B(u)\cap S(w)=\emptyset$.\\

\textbf{Case 2.2 : $\mathbf{w\in A(v)}$ :}\\
Suppose that $v$ gets attached to corner $q\in K(u)$ and edge $E_l(q)$. (See Section~\ref{Comments}). Clearly $q =
p(u)+r(u).S_1 = r(u).S_1$ for some $S_1\in {\cal S}$. Suppose that $w$ gets attached to corner $q'\in K(v)$ and edge
$E_{l'}(q')$ where $l'= J(v)$. Clearly $q'=p(v) + r(v).S_2$ for some $S_2\in {\cal S}$. From Step 3.2.2 of the algorithm, we
know that $q'\in K(v)\cap B(u)$. By Lemma~\ref{Two Corners}, the 2 corners of $B(v)$ which are inside $B(u)$ are given by
$p(v) + S'.r(v)$ and $p(v) + S'' r(v)$ where $S'$ and $S''$ are such that $S' = -S_1$ and $S''$ differs from $S_1$ in all
co-ordinates other than the $l^{th}$ co-ordinate. Thus the two strings $S_1$ and $S_2$ are related as follows:
$$S_2[m]\neq S_1[m], \forall m \neq l \mbox{ and } S_2[l] \mbox{ may or may not be equal to } S_1[l]$$
Let $|A(u)|=t$ and $|A(v)|=t'$. By Lemma~\ref{Co-ordinates}, $p(v) = r(u).S_1 - \frac{r(u)}{2}(1+ \frac{i}{t+1}).S_1[l].e_l +
S_1[j].\frac{r(v)}{2}.e_j$ for some $i \in [t]$ and $j = J(v) \neq l$. Again by Lemma~\ref{Co-ordinates}, $p(w) =
p(v)+r(v).S_2 - \frac{r(v)}{2}(1+ \frac{i'}{t'+1}).S_2[l'].e_{l'} + S_2[j'].\frac{r(w)}{2}.e_{j'}$ for some $i' \in [t']$ and
$j' = J(w) \neq l'$. Replacing value of $p(v)$ in previous equation, we get
\begin{eqnarray*}
p(w) &=& r(u).S_1 - \frac{r(u)}{2}(1+ \frac{i}{t+1}).S_1[l].e_l + S_1[j].\frac{r(v)}{2}.e_j \noindent\\
&+& r(v).S_2 - \frac{r(v)}{2}(1+ \frac{i'}{t'+1}).S_2[l'].e_{l'} + S_2[j'].\frac{r(w)}{2}.e_{j'} \noindent
\end{eqnarray*}
Also, recall from Step 3.2.2 of the algorithm that $l' =  J(v) = j$. Now, recalling that $j\neq l$, and thus $S_1[j] \neq
S_2[j]$ and also $j = l'$, we consider the distance between $p(u)$ and $p(w)$ along the $j^{th}$ co-ordinate: we get
$\rho(p(u),p(w)) = \rho(\mathbf{\overline 0}, p(w)) \geq |r(u) + \frac{r(v)}{2} - r(v) + \frac{r(v)}{2}(1+ \frac{i'}{t'+1})| =
r(u) + \frac{r(v)}{2}(\frac{i'}{t'+1})$. Now recall that $R(w) = 2r(w) = 2 \frac{r(v)}{8(t'+1)} < \frac{r(v)}{2(t'+1)} \leq
\frac{r(v)}{2}
\begin{bmatrix}
\frac{i'}{t'+1}
\end{bmatrix}$ for any $i' \in [t]$. Substituting this we get $\rho(p(u), p(w)) > r(u) + R(w)$ and we infer that $B(u)\cap S(w)=\emptyset$.
\hfill\qed
\end{proof}

\begin{lemma}
 If $uv\notin E(T)$, then $B(u) \cap B(v)= \emptyset$
\label{No Edge}
\end{lemma}
\begin{proof}
First suppose that none of $u$,$v$ is an ancestor of the other . Let $z$ be their least common ancestor. Let $z_u,z_v$ be
children of $z$ such that $u$ is a descendant of $z_u$ and $v$ is a descendant of $z_v$. By Lemma~\ref{Every descendant},
$B(u) \subseteq S(z_u)$ and $B(v) \subseteq S(z_v)$. But by Lemma~\ref{Superball of children}, $S(z_u) \cap S(z_v) = \emptyset$.
So, $B(u) \cap B(v) = \emptyset$.\\
Otherwise, without loss of generality, let $u$ be an ancestor of $v$. Since $uv\notin E(T)$ we know that $v$ is not a child of $u$.
Consider the path from $u$ to $v$ in $T$. Let $u,w,y$ be the first three vertices of this path. Clearly $v$ is a descendant of $y$ and therefore by
 Lemma~\ref{Every descendant}, $B(v) \subseteq S(y)$. Lemma~\ref{Grandchild} implies $S(y) \cap B(u) = \emptyset$.
 Therefore $B(v) \cap B(u) = \emptyset$
 \hfill\qed
\end{proof}

\begin{lemma}
$T$ is the intersection graph of the family $\{ B(u) | u\in V(T) \}$
\label{Intersection Graph}
\end{lemma}
\begin{proof}
Immediate from Lemma~\ref{Edge} and Lemma~\ref{No Edge}.
\hfill\qed
\end{proof}

\subsection{$T$ is the SIG of the family $\{ p(u) | u\in V(T) \}$}
\label{Step Two}

Consider the set ${\cal P}$ in $\mathbb{R}^d$ given by ${\cal P} = \{ p(u)\ |\ u\in V(T)\}$ computed by the algorithm. For
each $u\in V(T)$, let $r_u$ be distance of $p(u)$ from its nearest neighbor(s) in ${\cal P}$.

\begin{lemma}
 Let $u\in V(T)$ be a non-leaf vertex. Then there is a child $v$ of $u$ such that $\rho(p(u),p(v)) = r(u)$
\label{Non-leaves are good}
\end{lemma}
\begin{proof}
If $u$ has a leaf-child then we are through since any leaf-child of $u$ is placed at a corner of $B(u)$. If $u$ does not have
a leaf-child, then recall that we had designated a special child $v$ of $u$ as a pseudo-leaf. This pseudo-leaf is placed at a
corner of $B(u)$ in Step 3.2.1 of algorithm. So, $\rho(p(u),p(v)) = r(u)$
\hfill\qed
\end{proof}

\begin{lemma}
 Let $u\in V(T)$. Then there is a vertex $w\in V(T)$ such that $\rho(p(u),p(w)) = r(u)$
\label{All are good}
\end{lemma}
\begin{proof}
 If $u$ is not a leaf, then we are through by Lemma~\ref{Non-leaves are good}. If $u$ is a leaf, consider its parent $w$.
 Step 3.2.1 of the algorithm places $u$ at a corner of $B(w)$. So, $\rho(p(u),p(w)) = r(w) = r(u)$
\hfill\qed
\end{proof}

\begin{lemma}
 Let $uv\in E(T)$. Then $\rho(u,v)\geq max \{r(u),r(v)\}$.
\label{Maximum}
\end{lemma}
\begin{proof}
 Without loss of generality, let $u$ be the parent of $v$. Clearly, $r(u) \geq r(v)$. If $v\in C(u) - A(u)$, then
 $\rho(p(u),p(v)) = r(u)$. If $v\in A(u)$, then by Lemma~\ref{Co-ordinates}, $\rho(p(u),p(v)) = r(u) + \frac{r(v)}{2} > r(u)$.
 \hfill\qed
\end{proof}

\begin{lemma}
 Let $u\in V(T)$ . Then $r_u = r(u)$.
\label{Equal Radii}
\end{lemma}
\begin{proof}
Let $v\in V(T)$ be such that $uv\notin E(T)$. Then, by Lemma~\ref{Intersection Graph}, we have $\rho(p(u),p(v)) \geq r(u) +
r(v) > r(u)$.
 Let $w\in V(T)$ such that $uw\in E(T)$. Then, by Lemma~\ref{Maximum}, we have $\rho(p(u),p(v)) \geq r(u)$. Also, Lemma~\ref{All are good}
 guarantees us a vertex $y\in V(T)$ such that $\rho(p(u),p(y)) = r(u)$. Therefore, the  distance of $p(u)$ from its nearest neighbor(s)
 in ${\cal P}$ is $r(u)$.
 \hfill\qed
\end{proof}

\begin{lemma}
 $T$ is the SIG of the family $\{ p(u) | u\in V(T) \}$
\label{SIG Graph}
\end{lemma}
\begin{proof}
 Immediate from Lemma~\ref{Intersection Graph} and Lemma~\ref{Equal Radii}.
 \hfill\qed
\end{proof}

\subsection{SIG dimension of trees under $L_{\infty}$ metric}
\label{Theorem}

In the preceding sections we have shown that the set ${\cal P} = \{ p(u)\ |\ u\in V(T')\}$ given by the algorithm gives a
$d$-dimensional $SIG$ representation of $T$ for $d = \lceil \log_2(\beta + 2)\rceil$. Thus , we have $SIG(T) \leq \lceil
\log_2(\beta + 2)\rceil$ where $\beta = \beta(T)$. Recall that by Lemma~\ref{Lower Lemma}, $SIG(T) \geq \lceil \log (\beta +
1) \rceil$. We note that $\lceil \log_2(\beta + 1)\rceil = \lceil \log_2(\beta + 2)\rceil$ except when $\beta$ is one less
than a power of 2. Therefore we have the following theorem:\\

\hspace{-5mm}\textbf{\large{Theorem 1} :}\\
For any tree $T$, $SIG_\infty(T) \leq \lceil \log_2(\beta + 2)\rceil$ where $\beta = \beta (T)$.
If $\beta$ is not of the form $2^k-1$, for some integer $k \geq 1$, we have $SIG_\infty(T) = \lceil \log_2 (\beta + 2) \rceil$.

\section{When $\beta = 2^k - 1$ for some $d \geq 1$}
\label{Special Case}

\begin{figure}[!h]
\centering
\includegraphics[width=5.5in]{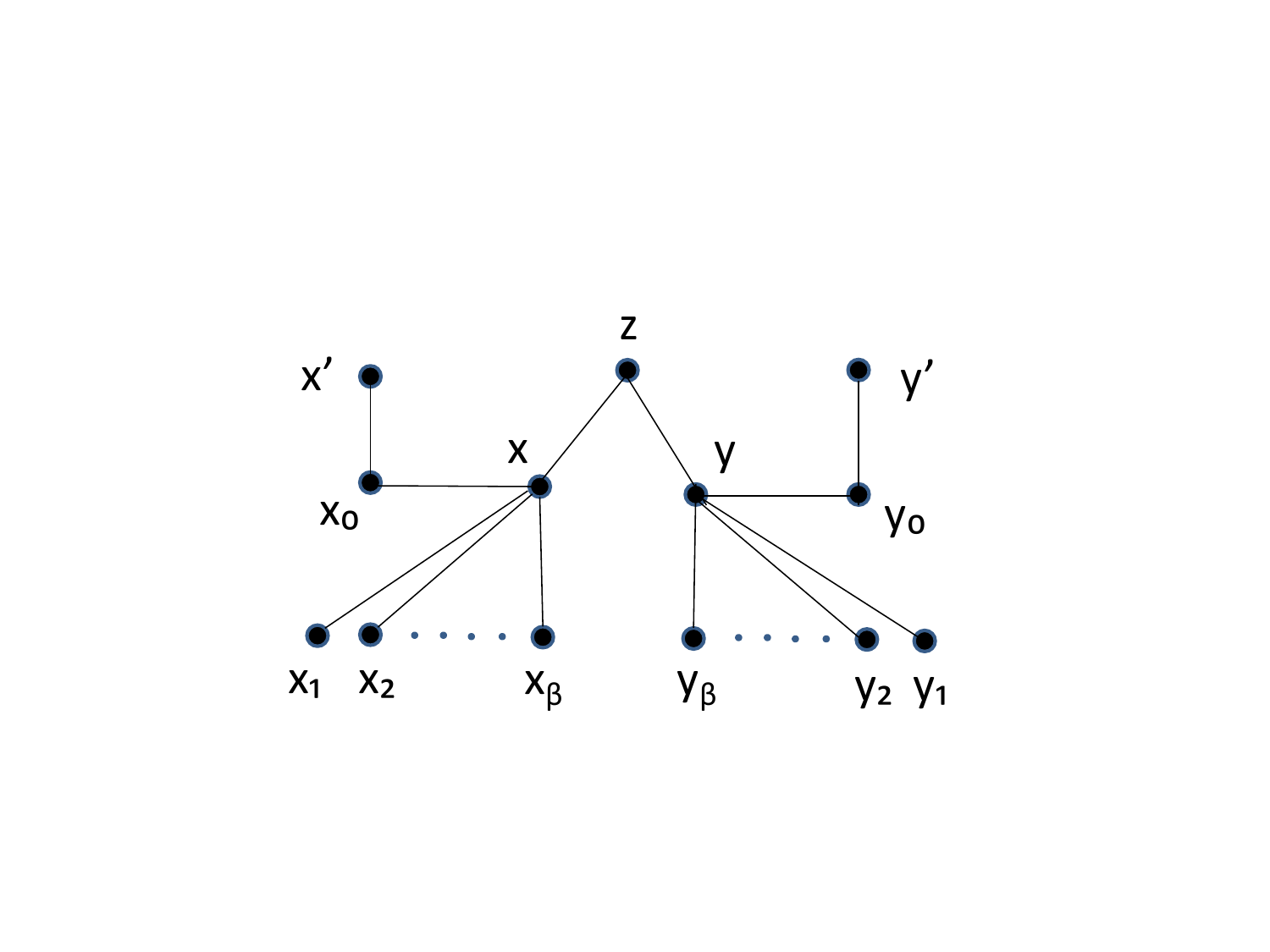}
\vspace{-15mm}
\caption{The graph $H$ with $\beta(H)=2^k-1$ and $SIG_\infty(H)=k+1$
\label{fig-tight-example}
}
\end{figure}

By Theorem 1 and Lemma~\ref{Lower Lemma}, we know that $(k+1) = \lceil \log_2(\beta+2)\rceil \geq SIG_{\infty}(T)
\geq \lceil \log_2(\beta+1)\rceil = k$. In this section we show that both values namely $k$ and $k+1$ are achievable.\\
~\\
\noindent{\bf Example where $\mathbf{SIG_\infty(T) = k}$ with $\mathbf{\beta(T)=2^k-1}$}\\
Consider a star graph on $b+2$ vertices, where
$b=2^k-1$. Since this is the complete bipartite graph with one vertex on one part and
$b+1$ vertices on the other part we denote it as $K_{1,b+1}$. Recalling the definition of $\beta(T)$ from
Section~\ref{Leaf-degree}, we note that $\beta(K_{1,b +1}) = \beta= b=2^k-1$. From Theorem 7 of~\cite{Michael}, we know that
$SIG_{\infty}(K_{1,b +1}) = \lceil \log_2(b +1)\rceil = k$. \\
~\\
\noindent{\bf Example where $\mathbf{SIG_\infty(T) = k + 1}$ with $\mathbf{\beta(T)=2^k-1}$}\\
Consider the tree $H$ illustrated in Fig.~\ref{fig-tight-example} on $(2b + 7)$ vertices, where $b=2^k-1$.
Recalling the definition of $\alpha$ and $\beta$ from Section~\ref{Leaf-degree} , we note that $\alpha(H)=\beta(H)=\beta=b=2^k-1$.
we show that its $SIG$ dimension
is $k+1$.

~\\


\vspace{-12mm}

\begin{lemma}
$SIG_{\infty}(H) = (k+1)$
\label{d+1}
\end{lemma}
\begin{proof}

By
Theorem 1, we have $SIG_{\infty}(H) \leq \lceil \log(2^k-1+2) \rceil=k+1$. Now we will show that $SIG_\infty(H)\geq k+1$. In a
$SIG$ graph, every vertex is adjacent to all of its nearest neighbor. Since the only two neighbors of $z$ are $x$ and $y$
either $x$ or $y$ or both must be nearest neighbor(s) of $z$. Without loss of generality, let $x$ be the nearest neighbor of
$z$. Also, for each $x_i$, $1 \leq i \leq \beta$, $x$ is the nearest neighbor, since $x$ is the only vertex adjacent to $x_i$.
Note that $\{x_0,z,x_1,x_2,\ldots,x_\beta\}$ form an independent set and $x$ is the nearest neighbor of each vertex from the
set $\{z,x_1,x_2,\ldots,x_\beta\}$. Let $SIG_{\infty}(H) = t$. Doing a similar analysis as in the proof of Lemma~\ref{Lower
Lemma}, we get $2^t > (\beta + 1)$ i.e  $2^t > 2^k$. Therefore, $t \geq (k + 1)$. It follows that $t = SIG_\infty(H) = k + 1$.
\hfill\qed
\end{proof}

~\\

\end{document}